\documentclass[11pt]{article}
\usepackage{graphicx}
\usepackage{hyperref}
\usepackage{color}
\usepackage{amssymb,latexsym,amsmath}
\usepackage{epstopdf}
\newtheorem{lemma}{Lemma}[section]
\newtheorem{theorem}[lemma]{Theorem}
\newtheorem{prop}[lemma]{Proposition}
\newtheorem{claim}[lemma]{Claim}
\newtheorem{cor}[lemma]{Corollary}
\newtheorem{remark}[lemma]{Remark}

\newtheorem{quest}[lemma]{Question}

\newcommand{\pf}{\noindent{\em Proof: }}
\newcommand{\proof}{\noindent{\em Proof: }}
\newcommand{\epf}{\hfill\hbox{\rule{3pt}{6pt}}\\}
\newcommand{\forme}[1]{}
\newcommand{\HH}{{\texttt h}}

\begin{document}

\date{\today}
\title{Distance-regular graphs without $4$-claws}
\author{{\bf Sejeong~Bang }\\Department of Mathematics\\ Yeungnam University, Gyeongsan-si, 
Gyeongbuk 38541, Republic of Korea \\
e-mail: sjbang@ynu.ac.kr\\
\\
{\bf Alexander Gavrilyuk}\\ School of Mathematical Sciences\\ 
University of Science and Technology of China, Hefei, 230026, Anhui, PR China;\\
N.N. Krasovskii Institute of Mathematics and Mechanics,\\ 
Ural Branch of Russian Academy of Sciences, Yekaterinburg, Russia\\
e-mail: alexander.gavriliouk@gmail.com\\
\\
{\bf Jack Koolen}\\ School of Mathematical Sciences\\ 
University of Science and Technology of China, Hefei, 230026, Anhui, PR China \\
e-mail: koolen@ustc.edu.cn}
\maketitle


\begin{abstract}
We determine the distance-regular graphs with diameter at least $3$ and $c_2\geq 2$ 
but without induced $K_{1,4}$-subgraphs.
\end{abstract}

\section{Introduction}

A $t$-{\it claw} is a complete bipartite graph $K_{1,t}$ with parts of size $1$ and $t$. 
We say that a graph is $t$-{\it claw-free} 
(or {\it without $t$-claws}), if it does not contain a $t$-claw as an induced subgraph.
For further definitions and notations, we refer the reader to Section 2.

Much attention has been paid to the study of $3$-claw-free graphs as a generalization 
of line graphs (it is easy to see that line graphs are $3$-claw-free).
The most general result in this direction belongs to Chudnovsky and Seymour \cite{ChudnovskySeymour},
who showed in a series of papers that every connected $3$-claw-free graph 
can be obtained from one of the so-called basic $3$-claw-free graphs by some simple operations.
Kabanov and Makhnev \cite{KabanovMakhnev} and (independently) Blokhuis and Brouwer \cite{no3-claws} 
determined all distance-regular graphs without 3-claws.

It is a well-known fact that the smallest eigenvalue of a line graph is at least $-2$.
In fact, any connected graph has smallest eigenvalue at most $-1$ with equality if and only 
if the graph is complete, i.e., its vertices are pairwise adjacent. 
It was shown by Cameron et al. \cite{CameronEtAl} (cf. \cite[Theorem~3.12.2]{bcn}) that 
a connected regular graph with smallest eigenvalue at least $-2$ is either a line graph or 
a cocktail party graph, or the number of its vertices is at most 28.

We recall that a {\it clique} in a graph is a set of vertices that induces a complete subgraph.
If a graph $\Gamma$ is regular with valency $k$, then its line graph ${\cal L}(\Gamma)$ 
is regular with valency $2(k-1)$, and it contains a set ${\cal C}$ of maximal cliques of 
the same order equal to $1+\frac{2(k-1)}{2}$ such that every edge of ${\cal L}(\Gamma)$ 
belongs to a unique member of ${\cal C}$. 

For distance-regular graphs, this observation can be formulated in a more general way.
Let $\Gamma$ be a distance-regular graph with valency $k$, diameter $D\geq 2$ and 
smallest eigenvalue $\theta_D$. 
Then any clique $C$ of $\Gamma$ satisfies the inequality
\begin{equation}\label{dc-bd}
|C|\leq 1+\frac{k}{|\theta_{D}|}, 
\end{equation}
which is known as the {\it Delsarte bound} (see \cite[Proposition 4.4.6]{bcn}).
This was shown by Delsarte \cite{Delsarte} for strongly regular graphs, and Godsil \cite{godsil-93-paper}
generalised it to distance-regular graphs.
A clique $C$ in $\Gamma$ is called a {\it Delsarte clique}, if $C$ contains exactly 
$1+\frac{k}{|\theta_D|}$ vertices.
Note that, if $\Gamma$ contains a Delsarte clique, 
then $\theta_D$ must be integral as $1+\frac{k}{|\theta_D|}$ is integral.

Godsil \cite{godsil-93-paper} introduced the following notion of a geometric 
distance-regular graph: a non-complete distance-regular graph $\Gamma$ is called {\it geometric}, 
if there exists a set of Delsarte cliques ${\cal C}$ such that every edge of $\Gamma$ lies 
in a unique member of ${\cal C}$. In this case we say that $\Gamma$ is {\it geometric 
with respect to} $\mathcal{C}$. In case of diameter 2, this notion is equivalent to 
that of geometric strongly regular graphs as introduced by Bose \cite{Bose}.
Examples of geometric distance-regular graphs include the Hamming graphs, the Johnson graphs,
the Grassmann graphs, the dual polar graphs, the bilinear forms graphs.

Note that the notion of a geometric distance-regular graph 
is stronger than those of a {\it clique geometry}, as defined by Metsch \cite{Metsch}, 
or an {\it asymptotic Delsarte geometry}, studied by Babai and Wilmes \cite{Babai}. Both definitions 
played an important role in recent progress on the complexity of the graph isomorphism problem 
(restricted to the class of strongly regular graphs) \cite{BabaiSRG}, and on the problem of classifying 
primitive coherent configurations with large automorphism groups \cite{Wilmes}.
In particular, Spielman \cite{Spielman} improved the complexity of isomorphism testing of 
strongly regular graphs, found by Babai \cite{Babai80}, using the following result by 
Neumaier \cite{neumaier-m}: for a fixed integer $\theta\geq 2$, 
there are only finitely many coconnected non-geometric distance-regular graphs with smallest eigenvalue 
at least $-\theta$ and diameter $2$. Koolen and Bang \cite{-m} generalised this as follows: 
there are only finitely many coconnected non-geometric distance-regular graphs with 
smallest eigenvalue at least $-\theta$ and with given diameter $D\geq 2$ or having the intersection 
number $c_2\geq 2$ (in fact, the Bannai-Ito conjecture recently proved by Bang et al. \cite{BIConjecture} 
shows that the latter condition on $c_2$ is not necessary). 

Neumaier \cite{neumaier-m} also showed that except for a finite number of graphs, 
all coconnected strongly regular graphs with a given smallest eigenvalue are geometric, 
and they are either Latin square graphs or Steiner graphs. 
In \cite{Wilson}, Wilson showed that there are super-exponentially many Steiner 
graphs, and similarly there are super-exponentially many Latin square graphs for certain parameter sets, 
see. This shows that the above-mentioned result of Neumaier is the best we can hope 
for the case of distance-regular graphs of diameter two. The situation for distance-regular graphs 
with larger diameter seems to be different: Koolen and Bang \cite{-m} conjectured that, 
for a fixed integer $\theta\geq 2$, any geometric distance-regular graph with smallest
eigenvalue $-\theta$, diameter $D\geq 3$ and $c_2\geq 2$ is either a Johnson graph, a Grassmann graph, 
a Hamming graph, a bilinear forms graph or the number of vertices is bounded above by a function
of $\theta$.

Note that a geometric distance-regular graph with smallest  eigenvalue $-t$ is $(t+1)$-claw-free, 
since every edge is contained in exactly one Delsarte clique, 
and hence each vertex is contained in exactly $t$ Delsarte cliques.

On the other hand, Godsil \cite{godsil-93-paper} showed that $\Gamma$ is geometric, 
if $\Gamma$ is $t$-claw-free for any $t>|\theta_D|$ and 
the intersection numbers $a_1$ and $c_2$ of $\Gamma$ satisfy the following inequality:
\[
{a_1+1>(-2\theta_D-1)(c_2-1).}
\]

Koolen and Bang \cite{-m} showed that $\Gamma$ is geometric, if $a_1\geq \lfloor \theta_D\rfloor^2c_2$ holds. 
Moreover, they \cite[Proposition~9.7]{drgsurvey} showed that a distance-regular graph $\Gamma$ of 
diameter $D\geq 2$ having valency $k\geq {\rm max}\{t^2-t,\frac{t^2-1}{t}(a_1+1)\}$ for some 
integer $t\geq 3$, is geometric with smallest eigenvalue $-t$ if and only if $\Gamma$ is $(t+1)$-claw-free.

Clearly, a connected graph without 2-claws is a complete graph. 
As we mentioned above, distance-regular graphs without 3-claws were classified in 
\cite{KabanovMakhnev} and \cite{no3-claws}.

In this paper, we classify distance-regular graphs without $4$-claws. 
Hiraki, Nomura and Suzuki \cite{k6} and Yamazaki \cite{yamazaki} considered distance-regular graphs 
that are locally a disjoint union of three cliques of size $a_1+1$ (i.e., 4-claw-free), 
and these graphs for $a_1\geq 1$ turn out to be geometric with smallest eigenvalue $-3$ 
(see Remark \ref{geo-rmk}). The problem of classification of such graphs seems to be very hard. 
Hence we may assume $c_2\geq 2$.
Bang \cite{gdrg-3} showed that if $k>\frac{8}{3}(a_1+1)$, then a distance-regular graph of valency $k$ 
is $4$-claw-free if and only if $\Gamma$ is geometric with smallest eigenvalue $-3$.
In \cite{gdrg-3}, Bang started the classification of geometric distance-regular graphs with 
smallest eigenvalue $-3$ and $c_2>1$, and the results of Bang and Koolen \cite{non-exist} 
and Gavrilyuk and Makhnev \cite{45} completed the classification. 

Our main result is the following classification of $4$-claw-free distance-regular graphs.

\begin{theorem}\label{theo-main}
Let $\Gamma$ be a distance-regular graph with diameter $D\geq 3$, valency $k>3$ but without $4$-claws. 
If $c_2\geq 2$, then $\Gamma$ is one of the following graphs:
\begin{itemize}
\item[(1)] The Taylor graph with $\iota(\Gamma)=\{k,\frac{k-1}{2},1;1,\frac{k-1}{2},k\}$ with 
$k\in \{5,13,17\}$;
\item[(2)] The Klein graph with $\iota(\Gamma)=\{7,4,1;1,2,7\}$;
\item[(3)] The Gosset graph with $\iota(\Gamma)=\{27,10,1;1,10,27\}$;
\item[(4)] The halved $e$-cube with $\iota(\Gamma)=\{{e\choose 2},{e-2\choose 2},1;1,6,15\}$ with $e\in \{6,7\}$;
\item[(5)] The Hamming graph $H(3,q)$ with $q\geq 3$;
\item[(6)] The Johnson graph $J(n,3)$ with $ n \geq 6$;
\end{itemize}
or $\Gamma$ is a putative distance-regular graph with one of the following intersection arrays:
\begin{itemize}
\item[(7)] $\iota(\Gamma)=\{44,24,1;1,12,44\}$;
\item[(8)] $\iota(\Gamma)=\{64,34,1;1,17,64\}$;
\item[(9)] $\iota(\Gamma)=\{104,54,1;1,27,104\}$.
\end{itemize}
\end{theorem}

The proof of Theorem \ref{theo-main} is organised 
according to the following three cases with respect to the valency $k$ 
and the intersection number $a_1$ of a distance-regular graph $\Gamma$ 
satisfying the condition of the theorem:
\begin{itemize}
\item[(T)] $k\leq 2(a_1+1)$,
\item[(N)] $2(a_1+1)<k\leq \frac{8}{3}(a_1+1)$,
\item[(G)] $k>\frac{8}{3}(a_1+1)$.
\end{itemize}

If $\Gamma$ satisfies case (T), then it is either 
the Johnson graph $J(7,3)$ or the halved $7$-cube 
or a Taylor graph (see Proposition \ref{prop-a1-half-valency}). 
In Section \ref{sect-taylor}, we classify the Taylor graphs without $4$-claws.

Now we consider the distance-regular graphs that satisfy cases (G) or (N). 
Here we show that the diameter $D$ of $\Gamma$ is at most 4 
(see Lemma \ref{vka}), and the smallest eigenvalue $\theta_D$ of $\Gamma$ 
is bounded below in terms of the ratio $\frac{k}{a_1+1}$ (see Lemma \ref{terw-ev-lemma}).

If $\Gamma$ satisfies case (G), then 
it was essentially shown by Bang (see Proposition \ref{prop-large-valency})
that $\Gamma$ is $4$-claw-free if and only if it is geometric 
with smallest eigenvalue $-3$. 
If $\Gamma$ is geometric with case (N), then we observe that 
there are only finitely many possibilities for $\theta_D$, 
and this enables us to classify all geometric 
distance-regular graphs with $c_2\geq 2$ but without $4$-claws in Section \ref{sect-geometric}.

In Section \ref{sect-nongeom}, we consider non-geometric distance-regular graphs with case (N).
Since the smallest eigenvalue $\theta_D$ of $\Gamma$ 
is bounded below (see Lemma \ref{terw-ev-lemma}), it follows by \cite{-m} 
that there are only finitely many non-geometric distance-regular graphs. 
In particular, they showed that the valency $k$ is bounded above as $O(m_1^2)$, 
where $m_1$ is the multiplicity of second largest eigenvalue.
However, their bound for $k$ is too large, 
which makes the search for feasible intersection arrays for $\Gamma$ almost intractable.
In Subsection \ref{sect-valencybounds}, we obtain a better valency bound of order $O(m_1)$ 
by using that any clique size gives a lower bound for $m_1$ (see Lemma \ref{m1ineq}), 
while the size of a maximum clique in its turn 
can be bounded below by using the 4-claw-free property (see Lemma \ref{terw-ev-lemma}).
In Subsection \ref{sect-comp}, we find all feasible intersection arrays, whose 
valencies satisfy the bounds of Subsection \ref{sect-valencybounds}.
For most of the arrays we found, we show in Subsection \ref{sect-nonexist}
that there are no corresponding distance-regular graphs without $4$-claws.
Moreover, in Proposition \ref{39}, we show that there does not exist 
a distance-regular graph with intersection array $\{39,24,1;1,4,39\}$.

Some of our results can be possibly generalised for distance-regular graphs without 
$t$-claws for any $t$, and Section \ref{sect-prelim} contains 
some basic definitions from the theory of distance-regular graphs and 
preliminary results on intersection numbers and eigenvalues of 
a distance-regular graph without $t$-claws.

Note that the complement of a $t$-claw-free strongly regular graph 
does not contain any $t$-clique. In general, it appears to be hard to show 
the non-existence of $t$-cliques in a strongly regular graph by using 
only its parameters (for example, the Shrikhande graph and the $(4\times 4)$-grid both 
have the same parameters, but the latter one contains a $4$-clique, while the former one 
does not; see further \cite{4clique} for the case $t=4$).
This suggests that the assumption $D\geq 3$ in Theorem \ref{theo-main} is crucial.

\section{Definitions and preliminaries}\label{sect-prelim}
\subsection{Definitions and notation}
All graphs considered in this paper are finite, undirected and simple (for more background information,
see \cite{bcn}). For a connected graph $\Gamma=(V(\Gamma),E(\Gamma))$, the {\it distance} $\partial(x,y)$
between any two vertices $x,y$ is the length of a shortest
path between $x$ and $y$ in $\Gamma$, and the {\it diameter} $D$ is the
maximum distance ranging over all pairs of vertices of $\Gamma$. 

For any vertex $x$ of $\Gamma$, let $\Gamma_i(x)$ be
the set of vertices in $\Gamma$ at distance precisely $i$ from $x$, where $0 \leq i\leq D$.
For a non-empty subset $S\subseteq V(\Gamma)$, $\langle S\rangle$ denotes the {\it induced subgraph} on $S$.
The {\em local graph} of a vertex $x$ is $\langle \Gamma_1(x) \rangle $. 
For a class $\mathcal{G}$ of graphs, 
a graph $\Gamma$ is called {\it locally} $\mathcal{G}$, if, for any vertex of $V(\Gamma)$, 
its local graph is isomorphic to a graph from $\mathcal{G}$. 
If all graphs from $\mathcal{G}$ are isomorphic to a graph $\Delta$, 
we say that $\Gamma$ is {\it locally} $\Delta$.
A graph $\Gamma$ is {\it regular} with valency $k$, if $\Gamma_1(x)$ 
contains precisely $k$ vertices for all $x\in V(\Gamma)$.

For a set of vertices $x_1,\ldots,x_n$ of $\Gamma$, 
let $\Gamma_1(x_1,\ldots,x_n)$ denote $\cap_{i=1}^n \Gamma_1(x_i)$.
Recall that a {\em clique} ({\em coclique}, resp.) is a set of pairwise 
adjacent (non-adjacent, resp.) vertices. 
A connected non-complete graph $\Gamma$ is called {\em Terwilliger}, 
if for any two vertices $x,y$ at distance $2$,
$\langle \Gamma_1(x,y)\rangle$ is a clique, whose order does not depend 
on the choice of $x$ and $y$.

A connected graph $\Gamma$ with diameter $D$ is called {\em distance-regular},~
if there exist integers $b_{i-1}$, $c_{i}$ $(1\leq i\leq D$) such that, for any two 
vertices $x,y\in V(\Gamma)$ with $\partial(x,y)=i$, there are precisely $c_i$ neighbors 
of $y$ in $\Gamma_{i-1}(x)$ and $b_i$ neighbors of $y$ in $\Gamma_{i+1}(x)$. 
Set $c_0=b_{D}=0$. The numbers $b_i, c_i$ and
$a_i:=b_0-b_i-c_i~(0 \leq i\leq D)$ are called the {\em
intersection numbers} of $\Gamma$. We observe $a_0=0$ and $c_1=1$. The array 
\begin{equation}\label{eq-array}
\iota(\Gamma):=\{b_0,b_1,\ldots,b_{D-1};c_1,c_2,\ldots,c_{D}\}
\end{equation}
is called the {\em intersection array} of $\Gamma$. 
In particular, $\Gamma$ is a regular graph with valency $k:=b_0$. 
We further note that, for $1\leq i\leq D$, the number $|\Gamma_i(x)|$ 
does not depend on the choice of a vertex $x$ of $\Gamma$. Define $k_i:=|\Gamma_i(x)|$
for any vertex $x$ and $i=0,1,\ldots,D$. Then we have $k_0=1,~k_1=k,~c_{i+1}k_{i+1}=b_{i} k_{i}~~(0 \leq i\leq D-1)$ and thus
\begin{equation} \label{ki}
k_i=\frac{b_1\cdots b_{i-1}}{c_2 \cdots c_i}k~~~(2\leq i\leq D),
\end{equation}
so that $v:=|V(\Gamma)|=k_0+k_1+\ldots+k_D$.


A distance-regular graph with diameter 2 is called a {\it strongly regular} graph. 
We say that a strongly regular graph $\Gamma$ has {\it parameters} $(v,k,\lambda,\mu)$, 
where $\lambda:=a_1$, and $\mu:=c_2$.

Let $\Gamma$ be a distance-regular graph of diameter $D\geq 2$.
The {\em adjacency matrix} $A:=A(\Gamma)$ is the $v\times v$-matrix with rows and columns indexed by $V(\Gamma)$, where the $(x,y)$-entry of $A$ is $1$, if $\partial(x,y)=1$, 
and $0$ otherwise. The eigenvalues of $\Gamma$ are the eigenvalues of $A$.
It is well-known that a distance-regular graph $\Gamma$ with diameter $D$ has
exactly $D+1$ distinct eigenvalues $k=\theta_0>\theta_1>\cdots >\theta_D$, which are
the eigenvalues of the following tridiagonal matrix (cf. \cite[p.128]{bcn}):
\begin{equation*}\label{mtx-L}
 L_1:= \left\lgroup
 \begin{tabular}{llllllll}
 $0$ & $k$\\
 $c_1$ & $a_1$ & $b_1$\\
 & $c_2$ & $a_2$ & $b_2$\\
 & & . & . & .\\
 & & & $c_i$ & $a_i$ & $b_i$\\
 & & & & . & . & .\\
 & & & & & $c_{D-1}$ & $a_{D-1}$ & $b_{D-1}$\\
& & & & &\makebox{\hspace{.324cm}} &
 $c_{D}$ & $a_{D}$
 \end{tabular}
 \right\rgroup.
\end{equation*}
The {\it standard sequence} $\{u_i(\theta)\mid 0 \leq i\leq D\}$
corresponding to an eigenvalue $\theta$ is a sequence satisfying the following recurrence relation:
\[
c_iu_{i-1}(\theta) + a_iu_i(\theta) + b_iu_{i+1}(\theta) = \theta u_i(\theta) ~~(1 \leq i \leq D),
\]
where $u_0(\theta) = 1$ and $u_1(\theta) = \frac{\theta}{k}$. Then the multiplicity of eigenvalue $\theta$ is given by
\begin{equation}\label{mi}
m (\theta)=\frac{v}{\sum_{i=0}^{D}k_iu_i^2(\theta)},
\end{equation}
which is known as {\em Biggs' formula} (cf. \cite[Theorem 21.4]{biggs}, \cite[Theorem 4.1.4]{bcn}).
Let $m_i$ denote the multiplicity of eigenvalue $\theta_i~(0 \leq i\leq D)$. Then
\begin{equation}\label{trace}
\sum_{i=0}^{D}m_i \theta_i^{\ell}= tr(A^{\ell}) = \mbox{the~number~of~closed~walks~of~length~$\ell$~in~$\Gamma$}
~~(\ell \geq 1),
\end{equation}
where $tr(A^{\ell})$ is the trace of matrix $A^{\ell}$ (cf. \cite[Lemma 2.5]{biggs}).\\

Let $\mathbb{Z}$ and $\mathbb{N}$ denote the set of integers and the set of positive integers, 
respectively.

\subsection{Basic results}

In this subsection, we collect several basic results on distance-regular graphs.
For the rest of Section \ref{sect-prelim}, let $\Gamma$ be a distance-regular graph of diameter $D\geq 2$, 
and with intersection array given by Eq. (\ref{eq-array}), valency $k$, and $D+1$ distinct 
eigenvalues $\theta_0=k$ $>$ $\theta_1$ $>$ $\ldots$ $>$ $\theta_D$ with multiplicities $m_0=1$, 
$m_1$, $\ldots$, $m_D$, respectively.

\begin{prop}(\cite[Theorem 16]{koo-park})\label{prop-a1-half-valency}
If $D\geq 3$ and $k\leq 2(a_1+1)$, then $\Gamma$ is one of the following graphs:
\begin{itemize}
\item[(1)] A polygon;
\item[(2)] The line graph of a Moore graph;
\item[(3)] The flag graph of a regular generalized $D$-gon of order $(s,s)$ for some $s$;
\item[(4)] A Taylor graph;
\item[(5)] The Johnson graph $J(7,3)$;
\item[(6)] The halved $7$-cube.
\end{itemize}
\end{prop}

\begin{lemma}\label{nonterwkbd}
If $k>2(a_1+1)$ and $\Gamma$ is not a Terwilliger graph, then 
$c_i-b_i\geq c_{i-1}-b_{i-1}+a_1+2$ holds (for all $1\leq i\leq D$) 
and 
\[k\geq 3 a_1-3c_2+7.\]
\end{lemma}
\pf
The inequality on $c_i-b_i$ follows from \cite[Theorem 5.2.1]{bcn}. 
Suppose that $\langle \{x,z_1,y,z_2\}\rangle$ is an induced quadrangle where $\partial(x,y)=\partial(z_1,z_2)=2$.
Put $X_i:=\{z_i\}\cup \Gamma_1(x,z_i)$ for $i=1,2$.
Then $2(a_1+1)-(c_2-2)\leq |X_1\cup X_2|\leq 2(a_1+1)$. As $k>2(a_1+1)$, there exists a vertex $z_3\in \Gamma_1(x)\setminus (X_1\cup X_2)$ and we find
\[ \left| \Gamma_1(x, z_3) \setminus (X_1\cup X_2)\right|\geq a_1-2(c_2-1). \]
Hence 
\[k=|\Gamma_1(x)|\geq |\{z_3\}|+|X_1\cup X_2|+|\Gamma_1(x, z_3) \setminus (X_1\cup X_2)|\geq 3a_1-3c_2+7,\]
which shows the result.\epf

The following results relate eigenvalues and intersection numbers of $\Gamma$.

\begin{theorem} (\cite[Theorem 3.6]{kyp})\label{kypeq}
The following inequality holds:
\[(\theta_1+1)(\theta_D+1)\leq -b_1,\]
with equality if and only if $D=2$.
\end{theorem}

\begin{lemma}(\cite[Proposition 3.2]{kyp})\label{KYP}
If $D=3$, then the following holds.
\begin{itemize}
\item[(1)] $\theta_1>\max\{ a_1+1-c_2, a_3 -b_2\}$.
\item[(2)] $\theta_1\geq \min\{a_3, \frac{a_1+\sqrt{a_1^2+4k}}{2}\}$.
\item[(3)] $\theta_2$ lies between $-1$ and $a_3-b_2$.
\end{itemize}
\end{lemma}

\begin{lemma}\label{d=3}
If $\theta_1<a_1$, then $D\leq 3$, and $D=3$ implies $b_1 < c_3$.
\end{lemma}
\proof (1) Suppose $D\geq 4$. For any two vertices $x,y$ with $\partial(x,y)=4$, the subgraph 
$\langle\Gamma_1(x)\cup \Gamma_1(y)\rangle$ is the
disjoint union of two connected regular graphs with valency $a_1$. By the interlacing result (see \cite[Theorem 3.3.1]{bcn})
we have $\theta_0\geq a_1$ and $\theta_1\ge a_1$, a contradiction.

(2) By (1) and Lemma \ref{KYP} (2), we find that $D=3$ and $a_3\leq \theta_1<a_1$. Hence the result follows by $b_1=k-2-(a_1-1)\leq k-2-a_3=c_3-2<c_3$. \epf

The following result gives some bounds on the multiplicity of the second largest eigenvalue of $\Gamma$.

\begin{lemma}\label{m1ineq} (cf. \cite[Theorem~3]{Powers}, \cite[Theorem 4.4.4]{bcn})
Suppose that $D\geq 3$ holds.
\begin{itemize}
\item[(1)] If $\Gamma$ has a clique of size $c$, then $$m_1\geq c.$$
\item[(2)] If $\theta_1$ is irrational with multiplicity $m_1<k$, then $$m_1\geq \frac{k}{2}.$$
\end{itemize}
\end{lemma}
\pf (1) It is straightforward by \cite[Theorem~3]{Powers}.

(2) By \cite[Theorem 4.4.4]{bcn}, $\theta_1+1$ and $\theta_D+1$ are 
conjugate and thus $m_1=m_D$. Moreover, each local graph of $\Gamma$ 
has eigenvalues $-1-\frac{b_1}{\theta_1+1}$ and $-1-\frac{b_1}{\theta_D+1}$ 
with multiplicities at least $k-m_1$ and $k-m_D$, respectively. This shows that 
$k \geq (k-m_1)+(k-m_D)=2k-2m_1$, and hence we find $m_1\geq \frac{1}{2}k$.
\epf

\subsection{Some structure theory for distance-regular graphs without $(t+1)$-claws}
In this subsection, we give some results on distance-regular graphs without $(t+1)$-claws ($t\geq 3$).

\begin{lemma}\label{terw-ev-lemma}
Suppose that $\Gamma$ contains a $t$-claw but no $(t+1)$-claws for some integer $t\geq 3$. 
With $\delta:=\frac{k}{a_1+1}$, the following holds.
\begin{itemize}
\item[(1)] $\Gamma$ contains a clique with at least $ \left(\frac{2\delta-t}{t\delta}\right)k+1$ vertices. 
In particular, $$\theta_{D}\geq  -\frac{t \delta}{2\delta-t}.$$
\item[(2)] $c_2\geq 1+\frac{2(t-\delta)}{t(t-1)}(a_1+1)$.
\end{itemize}
\end{lemma}
\pf By assumption, $\Gamma$ contains a $t$-claw, say $\langle \{x,y_1,\ldots,y_t\} \rangle$, 
where $\partial(x,y_i)=1~(1\leq i\leq t)$. 

(1) For each $1\leq i\leq t$, we define a set $$S_i:= \Gamma_1(x, y_{i}) \setminus
\bigcup_{j\in \{1,2,\ldots,t\}\setminus \{i\}} \Gamma_1(x, y_j).$$
We first show the following claim.
\begin{claim}\label{cl1}
There is an integer $1\leq j \leq t $ such that $|S_{j}| \geq \left(\frac{2\delta-t}{t\delta}\right)k-1$.
\end{claim}
\noindent{\em Proof of Claim \ref{cl1}: }
Let $\epsilon$ be the number of edges between $Y:=\{y_1,\ldots,y_t\}$ and $\Gamma_1(x)\setminus Y$. 
As each vertex in $Y$ has exactly $a_1$
neighbors in $\Gamma_1(x)$, we have $\epsilon =t a_1$.
Since each vertex in $S:=\cup_{i\in \{1,\ldots,t\}} S_i$ has a unique neighbor in $Y$
and each vertex in $\Gamma_1(x)\setminus (Y\cup S)$ has at least two neighbors in $Y$,
we find $\epsilon \geq |S|+2|\Gamma_1(x) \setminus (Y\cup S)|$.
Hence it follows by $a_1 = \frac{k}{\delta}-1$ and
$ |S|+|\Gamma_1(x) \setminus (Y\cup S)|=|\Gamma_1(x)\setminus Y|=k-t$ that
\[  2(k-t)-|S| = |S|+2|\Gamma_1(x) \setminus (Y\cup S)|\leq  \epsilon =t a_1 = t\left( \frac{k}{\delta}-1 \right) \]
and thus $|S|\geq \left(2-\frac{t}{\delta}\right)k-t$.
Now the claim follows as
$|S|=\sum_{i=1}^{t}|S_i|\geq \left(2-\frac{t}{\delta}\right)k-t$.
\epf

Without loss of generality, let $S_1$ satisfy $|S_1|\geq \left(\frac{2\delta-t}{t\delta}\right)k-1$ by Claim \ref{cl1}. If there exist two
non-adjacent vertices $v_1,v_2\in S_1$, then $\langle\{x, v_1,v_2, y_2,\ldots,y_t\}\rangle$ is a $(t+1)$-claw, which is impossible.
So $S_1\cup \{x,y_1\}$ induces a clique and thus
\[\left(\frac{2\delta-t}{t\delta}\right)k +1 \leq  |S_1\cup \{x,y_1\}|\leq 1+\frac{k}{|\theta_{D}|}\]
follows by Eq. (\ref{dc-bd}). Now we have 
\[\theta_{D}\geq  -\frac{t\delta}{2 \delta-t},\]
and this shows (1).

(2) Since there are no $(t+1)$-claws in $\Gamma$, we have 
$\Gamma_1(x)=\cup_{1\leq i\leq t}\left( \Gamma_1(x, y_i) \cup \{y_i\} \right)$.
It follows by the principle of inclusion and exclusion that
\[\delta (a_1+1)= k\geq \sum_{i=1}^{t}|\Gamma_1(x, y_i) \cup \{y_i\}|-\sum_{1\leq i<j\leq t}|
\Gamma_1(x,y_i,y_j)|\geq t(a_1+1)-\frac{t(t-1)}{2}(c_2-1).\]
This completes the proof.  \epf

In the following result, we find a slightly better bound for clique size when $t=3$.
Define $\mu_{\min}$ by 
$
\mu_{\min}:=\min \{|\Gamma_1(x,y,z)| \mid x\in V(\Gamma) \mbox{~and~} y,z\in \Gamma_1(x) \mbox{~satisfying~} \partial(y,z)=2\}.
$

\begin{cor}\label{taylorclaim}
Suppose that $\Gamma$ contains a $3$-claw but no $4$-claws. 
Then there exists a clique of size at least $k-1-2a_1+\mu_{\min}$.
\end{cor}
\pf
Let $\{x,y_1,y_2,y_3\}$ be a set, which induces a $3$-claw, where $y_i\in \Gamma_1(x),~i=1,2,3$. 
Put $\gamma:=|\Gamma_1(x,y_1,y_2,y_3)|$. As there are no $4$-claws, $\Gamma_1(x)=\cup_{i=1}^{3}(\Gamma_1(x,y_i)\cup \{y_i\})$ holds. Thus, $\gamma=k-3(a_1+1)+\sum_{1\leq i<j\leq 3}|\Gamma_1(x,y_i,y_j)|$ follows by 
\begin{eqnarray*}
k&=&|\Gamma_1(x)|\\&=&|\{y_1,y_2,y_3\}|+|\Gamma_1(x,y_1)|+|\Gamma_1(x,y_2)\setminus \Gamma_1(x,y_1)|+
|\Gamma_1(x,y_3)\setminus (\Gamma_1(x,y_1)\cup \Gamma_1(x,y_2))|\\&=&3+a_1+(a_1-|\Gamma_1(x,y_1,y_2)|)+(a_1-|\Gamma_1(x,y_2,y_3)|-|\Gamma_1(x,y_1,y_3)|+\gamma).
\end{eqnarray*}
Since there are no $4$-claws in $\Gamma$, we see that 
$X:=\Gamma_1(x,y_1)\setminus (\Gamma_1(x,y_2)\cup \Gamma_1(x,y_3))$ induces a clique, and hence $C:=\{x,y_1\}\cup X$ induces a clique of size $k-2a_1-1+|\Gamma_1(x,y_2,y_3)|$, which proves the claim.
\epf

\begin{lemma}\label{vka}
Suppose that $\Gamma$ contains a $3$-claw but no $4$-claws. 
If $D\geq 3$ and $3\leq k < 3(a_1+1)$, then the following holds.
\begin{itemize}
\item[(1)] $\Gamma$ contains an induced quadrangle (i.e., $c_2\geq 2$).
\item[(2)] $b_2\leq \min\{a_1+1, a_3+1, \frac{k}{-\theta_D} \}$.
\item[(3)] $D\leq 4$.
\item[(4)] If $D=3$, then $\theta_2\ge -1$.
\end{itemize}
\end{lemma}
\pf
(1) As $\Gamma$ contains a $3$-claw and $k<3(a_1+1)$, it follows from 
Lemma \ref{terw-ev-lemma}~(2) that $c_2\geq 2$. 
Suppose that $\Gamma$ has no induced quadrangles (i.e., $\Gamma$ is a  Terwilliger graph). 
Then, by \cite[Proposition 6 (3)]{koo-park}, $\Gamma$ is the icosahedron, the Conway-Smith graph, 
or the Doro graph. Both the Conway-Smith graph and the Doro graph have $4$-claws, and 
the icosahedron does not have a $3$-claw, a contradiction. 

(2) By (1) there is a quadrangle, 
say $\{x,z_1,y,z_2\}$, where $\partial(x,y)=2=\partial(z_1,z_2)$. As $D\geq 3$ and
$\Gamma$ is $4$-claw-free, $\langle \Gamma_3(x)\cap \Gamma_1(y)\rangle$ must be a clique of size $b_2$. 
This shows that $a_1\geq b_2-1$ and $a_3\geq b_2-1$. Since 
$\langle (\Gamma_3(x)\cap \Gamma_1(y))\cup \{y\}\rangle$ is a clique of size $b_2+1$, 
the result follows by Eq. (\ref{dc-bd}).

(3) Suppose $D\geq 5$. By (1), there is a quadrangle, say $\langle \{y,s,t,w \}\rangle$ with $\partial(y,t)=\partial(s,w)=2$.
As $D\geq 5$, there is a vertex $x$ satisfying $\partial(x,y)=3$ and $\partial(x,t)=5$. As $\{s,w \}$ is a coclique in $\Gamma_1(y)$ and $\Gamma$ is $4$-claw-free,
the set $C:=\Gamma_2(x)\cap \Gamma_1(y)$ induces a clique of size $c_3$. It follows by (2), 
Lemma \ref{nonterwkbd}, and $b_3\geq c_2$ that
$c_3\geq c_2-b_2+a_1+2+b_3\geq c_2+1+b_3\geq 2c_2+1$ and thus
\[c_3>2c_2.\]
Let $y_1\in C$ and take $x_1\in \Gamma_1(x, y_1)$. As $c_3>2c_2$, there is a vertex 
$y_2 \in C \setminus \Gamma_1(x_1)$,
and take a vertex $x_2\in \Gamma_1(x,y_2)$. Similarly, there is a vertex 
$y_3\in C \setminus(\Gamma_1(x_1)\cup \Gamma_1(x_2))$ as $c_3>2c_2$.
Take a vertex $x_3\in \Gamma_1(x,y_3)$. As the set $C$ induces a clique, for any $1\leq i< j\leq 3$, 
two vertices $x_i$ and $y_j$ has $c_2$ common neighbors in $C$ and thus $\{x_1,x_2,x_3\}$ is a coclique. On the other hand,
there exists a vertex $x_4\in \Gamma_1(x)\cap \Gamma_4(y)$ by $b_3\geq 1$, and the set $\{x_i\mid 1\leq i\leq 4\}$ induces a coclique as
$\partial(y,x_i)=2$ ($i=1,2,3$) and $\partial(y,x_4)=4$.
Hence $\langle \{x,x_i\mid 1\leq i\leq 4\}\rangle$ is a $4$-claw in $\Gamma$, which is impossible. This shows $D\leq 4$. 


(4) As $a_3-b_2\geq -1$ holds by (2), the result (4) follows from Lemma \ref{KYP} (3).   \epf

\begin{prop}\label{prop-large-valency} (cf. \cite[Theorem 3.2]{gdrg-3})
If $k>\frac{8}{3}(a_1+1)$, then the following statements are equivalent.
\begin{itemize}
\item[(1)] $\Gamma$ is $4$-claw-free.
\item[(2)] $\Gamma$ is a geometric distance-regular graph with smallest eigenvalue $-3$.
\end{itemize}
\end{prop}
\pf If $\frac{8}{3}(a_1+1)<k\leq 3$, then $k=3$ and $a_1=0$. 
Hence $\Gamma$ is the cube $H(3,2)$ (see \cite{k=3}), which satisfies both 
$(1)$ and $(2)$. If $k > \max\{3,\frac{8}{3}(a_1+1)\}$, then the result follows from \cite[Theorem 3.2]{gdrg-3}. \epf

\section{Taylor graphs without $4$-claws}\label{sect-taylor}
In this section we classify the Taylor graphs without $4$-claws (Theorem \ref{general-taylor}).\\

Let $\Gamma$ be a Taylor graph (i.e., a distance-regular graph with diameter three and intersection array $\iota(\Gamma)=\{k,c_2,1;1,c_2,k\}$).
The eigenvalues $\theta_i~(0\leq i\leq 3)$ of $\Gamma$ are
\begin{equation}\label{taylor-evs}
\theta_0=k,~~\theta_2=-1,~~\theta_i=\frac{k-2c_2-1+(-1)^{\frac{i-1}{2}}\sqrt{(k-2c_2-1)^2+4k}}{2}~~~(i=1,3),
\end{equation}
whose multiplicities $m_i=m(\theta_i)~(0\leq i\leq 3)$ are given by
\begin{equation}\label{taylorm1}
m_0=1,~m_1=\frac{k(k+1)}{\theta_1^2 +k},~m_2=k,~m_3=\frac{\theta_1^2 (k+1)}{\theta_1^2+k}.
\end{equation}

\begin{theorem}\label{general-taylor}
Let $\Gamma$ be a Taylor graph without $4$-claws. Then $\Gamma$ is one of the following graphs:
\begin{itemize}
\item[(1)] The hexagon; 
\item[(2)] The cube with $\iota(\Gamma)=\{3,2,1;1,2,3\}$;
\item[(3)] The icosahedron with $\iota(\Gamma)=\{5,2,1;1,2,5\}$;
\item[(4)] The Johnson graph $J(6,3)$ with $\iota(\Gamma)=\{9,4,1;1,4,9\}$;
\item[(5)] The Taylor graph with $\iota(\Gamma)=\{13,6,1;1,6,13\}$;
\item[(6)] The halved $6$-cube with $\iota(\Gamma)=\{15,6,1;1,6,15\}$;
\item[(7)] The Taylor graph with $\iota(\Gamma)=\{17,8,1;1,8,17\}$;
\item[(8)] The Gosset graph with $\iota(\Gamma)=\{27,10,1;1,10,27\}$.
\end{itemize}

Moreover, the hexagon, the cube and the Johnson graph $J(6,3)$ are the only geometric 
Taylor graphs without $4$-claws. 
\end{theorem}
\pf If $\Gamma$ is $3$-claw-free, then $\Gamma$ is either the hexagon or the icosahedron by \cite[Theorem 1.2]{no3-claws}. We assume that $\Gamma$ is a Taylor graph with $3$-claws but without $4$-claws. If $k=3$,  then $\Gamma$ is the cube (see \cite[Theorem 7.5.1]{bcn}). From now on, we assume $k\geq 4$. If $m_1=k$, then it follows by Eq. (\ref{trace}) and Eq. (\ref{taylorm1}) that $\theta_1=1$ and $\theta_3=-k$, which implies that $\Gamma$ is bipartite. This is impossible as $k\geq 4$ and $a_1=0$ (i.e., $\Gamma$ contains $4$-claws). Hence $m_1\neq k$. As $m_1+m_3=k+1$, $m_1 <k$ holds. For a vertex $x\in V(\Gamma)$, let
 $\Delta:=\langle \Gamma_1(x)\rangle$ be the local graph of $x$.
 By \cite[Theorem 1.5.3]{bcn}, $\Delta$ is a strongly regular graph with parameters $(\nu,\kappa,\lambda,\mu)=\left(k,a_1,a_1-\frac{1}{2}c_2-1,\frac{1}{2}a_1 \right)$ and eigenvalues $\eta_0=a_1,~\eta_1,~\eta_2$,  where $\eta_2<0<\eta_1$ and 
\begin{equation}\label{taylor2}
\eta_1+\eta_2=\frac{1}{2}a_1-\frac{1}{2}c_2-1 \mbox{~~and~~}
\eta_1 \eta_2=-\frac{1}{2}a_1.
\end{equation}
By Eqs. (\ref{dc-bd}), (\ref{taylor-evs}) and Corollary \ref{taylorclaim},
\begin{equation}\label{taylordb}
k-\frac{3}{2}a_1-1  \leq 1+\frac{k}{|\theta_3|}=1+\theta_1.
\end{equation}
If $m_1=m_3$, then $m_1=m_3=\frac{k+1}{2}$ and $\theta_1=-\theta_3=\sqrt{k}$ follow by Eqs. 
(\ref{taylor-evs}), (\ref{taylorm1}) and $tr(A)=0$
in Eq. (\ref{trace}). As $|V(\Gamma)|k a_1= tr(A^3)=\sum_{i=0}^{3}m_i\theta_i^3$ holds by Eq. (\ref{trace}), 
we find $a_1=c_2=\frac{k-1}{2}$ (thus $k$ is odd) and hence $k\leq 25$ holds as
\[ \frac{k-1}{4}=k-\frac{3}{2}a_1-1\leq 1+\theta_1=1+\sqrt{k} \]
by Eq. (\ref{taylordb}). 
On the other hand, $a_1$ is even and $k$ is a sum of two squares by $\mu=\frac12a_1 \in \mathbb{N}$ and \cite[Proposition 1.10.5 (i)]{bcn}. These restrictions imply $k\in \{5,9,13,17,25\}$. If $k=5$, then $\Gamma$ is the icosahedron, which is $3$-claw-free. If $k=9$, then $\Gamma$ is the Johnson graph $J(6,3)$ which contains $3$-claws but no $4$-claws. If $k=13,17$, then $\Delta$ is the Payley graphs $P(13)$ and $P(17)$, respectively. As the independence numbers of these Payley graphs are all three, $\Gamma$ contains $3$-claws but does not contain $4$-claws. 
If $k=25$, then each local graph of $\Gamma$ is a strongly regular graph with parameters $(25,12,5,6)$. Since all such stongly regular graphs contain $4$ cocliques (see \cite{srg25}), $\Gamma$ contains $4$-claws and thus $k\neq 25$. Now the results (4), (5) and (7) follow.\\ Now we assume that $m_1\neq m_3$. Then all the eigenvalues $\theta_i~(0\leq i\leq 3)$ are integral. As $m_1<k$ and by \cite[Theorem
4.4.4]{bcn}, $b_1=\ell (\theta_1+1)$ holds for $\ell\in \mathbb{N}$ with $1\leq \ell \leq
-\theta_3-1$, and
\begin{equation}\label{taylor1}
-1-\frac{b_1}{\theta_1+1}=-\ell -1
\end{equation}
is the smallest eigenvalue of $\Delta$, i.e., $\eta_2=-\ell-1$.
It follows by Eq. (\ref{taylor2}) with $\eta_2=-\ell-1$ that
\begin{equation}\label{taylor3}
c_2=\frac{\ell}{\ell +1}a_1+2 \ell \mbox{~~and~thus~~}k=(2\ell + 1) \left( \frac{a_1}{\ell +1}+1\right).
\end{equation}
By Eqs. (\ref{taylordb}), (\ref{taylor1}) and (\ref{taylor3}), we find that $(\ell -3)a_1 + 4(\ell^2 -1)\leq 0$, which shows $\ell \in \{1,2\}$.
If $\ell=1$, then $c_2=\frac{k+3}{3}$ and thus $\iota(\Gamma)=\{k,\frac{k+3}{3},1;1,\frac{k+3}{3},k\}$ with $\theta_1=\frac{k}{3}$ and $\theta_3=-3$ by Eq. (\ref{taylor-evs}). As $m_1=9-\frac{72}{k+9}\in \mathbb{N}$ holds by Eq. (\ref{taylorm1}), we find $k\in \{15,27,63\}$ as $m_1\neq m_3$.
By \cite[Corollary 1.15.3]{bcn}, $k\neq 63$. If $k=15$, then $\Gamma$ is the halved $6$-cube whose local graph is triangular graph $T(6)$ and independence number of $T(6)$ is three. If $k=27$, then $\Gamma$ is the Gosset graph whose local graph is the Schl\"{a}fli graph and its independence number is also three. 
If $\ell=2$, then it follows by Eq. (\ref{taylor3}) that $\iota(\Gamma)=\{k,\frac{2k+10}{5},1;1,\frac{2k+10}{5},k\}$ holds and thus $\theta_1=\frac{k}{5}$, $\theta_3=-5$ and $m_1=25-\frac{600}{k+25}$
by Eqs. (\ref{taylor-evs}) and (\ref{taylorm1}). As $m_1\in \mathbb{N}$, $m_1<k$, $m_1\neq m_3$ and $a_1$ is even, $k\in \{15,35,75,95,125,175,275,575\}$. Note that $\Delta ^c$ (the complement of local graph $\Delta$) is also a strongly regular graph with parameters $\left(k,\frac{2k+10}{5},\frac{k+25}{10},\frac{k+5}{5}\right)$. It follows by \cite[Theorem 3.3]{4clique} and \cite{4clique-web} that  for each  $k\in \{15,35,75,95,125,175,575\}$, $\Delta^c$ contains cliques of size $4$. This implies $\Delta$ contains cocliques of size $4$ and thus $\Gamma$ does contain $4$-claws, a contradiction. If $k=275$, then 
$\Delta$ is the McLaughlin graph which has a coclique of size $22$, this is also impossible. 
This completes the proof.
\epf

\section{Geometric distance-regular graphs without $4$-claws}\label{sect-geometric}
In this section we classify geometric distance-regular graphs with diameter $D\geq 3$ 
but without $4$-claws (see Theorem \ref{geo-thm}).\\

Let $\Gamma$ be a geometric distance-regular graph with respect to a set of Delsarte cliques $\mathcal{C}$ and diameter $D\geq 3$. For a vertex $x\in V(\Gamma)$ and a non-empty subset $S$ of $V(\Gamma)$, we denote $\partial(x,S)=\min\{\partial(x,y) \mid y\in S \}$. For a Delsarte clique $C\in \mathcal{C}$ and for each $i=0,1,\ldots,D-1$, let $ C_i:=\{x\in V(\Gamma) \mid \partial(x,C)=i\}$. For a vertex $x \in C_i$, define
\[\psi_i(x,C) := \left| \{ z \in C \mid \partial(x,z) = i\}\right|.\]
For $x, y\in V(\Gamma)$ with $\partial(x,y)=j~~(j=1,2,\ldots, D)$, define
\[ \tau_j(x,y;\mathcal{C}) :=\left| \{ C\in \mathcal{C} \mid x\in C \mbox{~~and~~} \partial(y,C)=j-1 \} \right|. \]
Since parameters $\psi_i(x,C)$ and $\tau_j(x,y;\mathcal{C})$ depend only on the distances $i$ and $j$ respectively (see \cite[Lemma 4.1]{DCG1} and \cite[Section 11.7]{godsil-93}), we denote $\psi_i:=\psi_{i}(x,C)~~(i=0,1,\ldots,D-1)$ and $\tau_j:=\tau_j(x,y;\mathcal{C})~(j=1,2,\ldots,D)$. 
We observe that for any vertex $x\in V(\Gamma)$ and for any Delsarte clique $C \in \mathcal{C}$, 
\begin{equation}\label{tau-ineq}
|\{X\in \mathcal{C} \mid x\in X \}|=-\theta_{D},~~\partial(x,C)\leq D-1~~ \mbox{~and~thus~} \tau_D=-\theta_{D}.
\end{equation}

The next lemmas summarize some results for geometric distance-regular graphs.

\begin{lemma}(\cite[Lemma 4.1, Lemma 4.2]{gdrg-3}, \cite[Lemma 5.1]{DCG1})\label{geo-basic-lem}
Let $\Gamma$ be a geometric distance-regular graph with diameter $D\geq 2$. Then the following holds.
\begin{itemize}
\item[(1)] $b_i = -(\theta_D+\tau_i)\left(1-\frac{k}{\theta_D}-\psi_i\right)~(1\leq i \leq D-1 )$.
\item[(2)] $c_i=\tau_i \psi_{i-1}~(1\leq i\leq D)$.
\item[(3)] $c_2=\psi_1$ if and only if $c_2=1$.
\item[(4)] If $D\geq 3$, then $\psi_1\leq \tau_2\leq -\theta_D-1$.
\end{itemize}
\end{lemma}

\begin{lemma}(\cite[Proposition 9.9]{drgsurvey})\label{surveyc2}
Let $t\geq 2$ be an integer, and let $\Gamma$ be a distance-regular graph with diameter $D\geq 2$ 
and valency $(t-1)(a_1+1)<k<t(a_1+t)$.
If $a_1>\frac{t(t+1)(c_2+1)}{2}$, then $\Gamma$ is geometric with $\theta_D=-t$.
\end{lemma}

\begin{lemma}\label{geo-d3}
Let $\Gamma$ be a geometric distance-regular graph with $3$-claws but without $4$-claws. 
Then either $k=3(a_1+1)$ with $c_2=1$ or $D\leq 3$ with $c_2\geq 2$.
\end{lemma}

\pf As $\Gamma$ is $4$-claw-free, $3\leq k\leq 3(a_1+1)$. It follows by Lemma \ref{vka} (1) that either $k=3(a_1+1)$ with $c_2=1$ or $3\leq k<3(a_1+1)$ with $c_2\geq 2$. Let $\Gamma$ satisfy $3\leq k<3(a_1+1)$ and $c_2\geq 2$. Then it is enough to show $D\leq 3$. We first show the following claim. 
\begin{claim} \label{geo-claim}
For any vertices $x$ and $y$ with $\partial(x,y)=2$, there is a quadrangle containing $x$ and $y$.
\end{claim}
\noindent{\em Proof of Claim \ref{geo-claim}:} Let $\partial(x,y)=2$. 
As $\Gamma$ is geometric, there is a Delsarte clique $C$ satisfying $x\in C$ and $\partial(y,C)=1$. 
As $c_2\geq 2$, we have $c_2>\psi_1$ by Lemma \ref{geo-basic-lem} (3) and thus there exists a vertex
$w \in \Gamma_1(x, y)\setminus C$. By $|\Gamma_1(w)\cap C|=|\Gamma_1(y)\cap C|=\psi_1$, $x\in \Gamma_1(w)\cap C $ and $x\not\in \Gamma_1(y)\cap C$, there exists a vertex $z \in \Gamma_1(y)\cap C$
such that $\partial(z,w)\ne 1$ and thus $\partial(z,w)=2$. Hence the induced subgraph on $\{x,z,y,w\}$ is a quadrangle. \epf

If $D\geq 4$, then there exist $3$ vertices $p,q,r\in V(\Gamma)$ satisfying $\partial(p,q)=\partial(q,r)=2$ and $\partial(p,r)=4$. By Claim \ref{geo-claim},
there are $4$ vertices $s_1,s_2\in \Gamma_1(p,q)$ and $s_3,s_4 \in \Gamma_1(q,r)$ such that $\langle \{p,s_1,q,s_2\}\rangle$ and
$\langle\{q,s_3,r,s_4\}\rangle$ are quadrangles, respectively. This is impossible as  the set $\{q,s_1,s_2,s_3,s_4\} $ induces a $4$-claw.
This shows $D\leq 3$.\epf

Now we prove the main result of this section.

\begin{theorem}\label{geo-thm}
Let $\Gamma$ be a geometric distance-regular graph with diameter $D\geq 3$, valency $k>3$ 
but without $4$-claws. Then $\Gamma$ is one of the following graphs:
\begin{itemize}
\item[(1)] One of the two generalized hexagons of order $(2, 2)$ with $\iota(\Gamma) = \{6,4,4;1,1,3\}$; 
\item[(2)] The generalized hexagon of order $(8,2)$ with $\iota(\Gamma) = \{24, 16, 16; 1, 1, 3\}$;
\item[(3)] The halved Foster graph with $\iota(\Gamma)=\{6,4,2,1;1,1,4,6\}$;
\item[(4)] A generalized octagon of order $(4,2)$ with $\iota(\Gamma)=\{12,8,8,8;1,1,1,3\}$;
\item[(5)] The Hamming graph $H(3,q)$ with $q\geq 3$;
\item[(6)] The Johnson graph $J(n,3)$ with $ n \geq 6$;
\item[(7)] The line graph of a Moore graph;
\item[(8)] The flag graph of a regular generalized $D$-gon of order $(s,s)$ for some $s$;
\item[(9)] A distance-regular graph with $\theta_D=-3$ satisfying $k=3(a_1+1)$, $c_2=1$ 
and $a_1,\frac{1}{3}c_D\geq 2$. 
\end{itemize}
\end{theorem}

\pf If $\Gamma$ is $3$-claw-free, then it follows by \cite[Theorem 1.2]{no3-claws} that $\Gamma$ is (7) or (8). (Note that these graphs are all geometric as they satisfy $k=2(a_1+1)$ and $\theta_D=-2$, and the icosahedron is non-geometric as it has irrational smallest eigenvalue.) Now we suppose that $\Gamma$ contains $3$-claws. If $k\leq 2(a_1+1)$, then it follows by Proposition \ref{prop-a1-half-valency} and Theorem \ref{general-taylor} that $\Gamma$ is the Johnson graph $J(n,3)$ with $n \in \{6,7\}$. If $k>\frac{8}{3}(a_1+1)$, then it follows by Proposition \ref{prop-large-valency}, \cite[Theorem 4.3]{gdrg-3} and \cite[Theorem 1.4]{non-exist} that $\Gamma$ is one of the following:\\
\noindent(a) One of the two generalized hexagons of order $(2, 2)$ with $\iota(\Gamma) = \{6, 4, 4; 1, 1, 3\}$; \\
\noindent(b) The generalized hexagon of order $(8,2)$ with $\iota(\Gamma)=\{24,16,16;1,1,3\}$;\\
\noindent(c) The halved Foster graph with $\iota(\Gamma)=\{6,4,2,1;1,1,4,6\}$;\\
\noindent(d) A generalized octagon of order $(4,2)$ with $\iota(\Gamma)=\{12,8,8,8;1,1,1,3\}$;\\
\noindent(e) The Hamming graph $H(3,q)$ with $q\geq 3$;\\
\noindent(f) The Johnson graph $J(n,3)$ with $n\geq 20$;\\
\noindent(g) A distance-regular graph with $\theta_D=-3$ satisfying $k=3(a_1+1)$, $c_2=1$ and $a_1,\frac{1}{3}c_D\geq 2$.\vskip0.001cm
From now on, we assume that $2(a_1+1)<k\leq \frac{8}{3}(a_1+1)$. Note that $D=3$ and $c_2\geq 2$ follow by Lemma \ref{geo-d3}. Since there is a quadrangle by Lemma \ref{vka} (1), there are two vertices $p$ and $q$ at distance two such that vertices
$p,q,r_1,r_2$ induce a quadrangle for some $r_1,r_2\in \Gamma_1(p,q)$. As $D=3$, there exists $r_3\in \Gamma_3(p)\cap \Gamma_1(q)$. This implies that
the set $\{q,r_1,r_2,r_3\}$ induces a $3$-claw and thus $\theta_3\in \{-3,-4,-5\}$ by Eq. (\ref{tau-ineq}) 
and Lemma \ref{terw-ev-lemma} (1). If $\psi_1=1$ then $k=(-\theta_3)(a_1+1)$, which is impossible as $k\leq \frac{8}{3}(a_1+1)$ and $\theta_3\leq -3$. Thus $\psi_1\geq 2$. If $\theta_3=-3$, then $\Gamma$ is the Johnson graph $J(n,3)$ with $8\leq n \leq 19$ by \cite[Theorem 4.3]{gdrg-3}. By Eq. (\ref{tau-ineq}) with 
$\theta_3\in \{-4,-5\}$, there are at least four Delsarte cliques
$C_{(i)}~(1\leq i\leq 4)$ containing a fixed vertex $x$.
If $|C_{(i)}|-1=\frac{k}{-\theta_3}\geq 3(\psi_1-1)+1$, then there exist vertices $y_1\in C_{(1)}$ and
$y_i\in C_{(i)}\setminus (\cup_{j=1}^{i-1} \Gamma_1(y_j))~(2\leq i\leq 4)$ such that the set 
$\{x,y_i\mid 1\leq i\leq 4\}$ induces a $4$-claw, which is impossible and hence
\begin{equation}\label{psiineq}
\frac{k}{-\theta_3}\leq 3(\psi_1-1).
\end{equation}
It follows by $a_1=(\frac{k}{-\theta_3}-1)+(-\theta_3-1)(\psi_1-1)$ with Eq. (\ref{psiineq}) 
that if $\theta_3=-4$, then
$$a_1=\frac{k}{4}-1+3(\psi_1-1)\geq \frac{k}{2}-1.$$ This is impossible as $k> 2(a_1+1)$. Hence $\theta_3=-5$ and thus $2\leq \psi_1\leq 4$ and $\psi_1\geq \frac{1}{15}k+1$ follow by Lemma \ref{geo-basic-lem} (4) and 
Eq. (\ref{psiineq}). Using Eq. (\ref{psiineq}), $\psi_1\in \{2,3,4\}$ and $2\psi_1\leq 1+\frac{k}{5}$ (see \cite[Lemma 5.2 (ii)]{DCG1}) we find $(\psi_1,k)\in \{(2,15),(3,25),(3,30),(4,35),(4,40),(4,45)\}.$ By the feasible conditions $k_2,k_3,b_1,b_2,c_2,c_3\in \mathbb{N}$, Eq. (\ref{trace}), Lemma \ref{geo-basic-lem} and $2\leq \psi_1<\psi_2$
(see \cite[Theorem 5.5]{DCG1}), we find that the feasible intersection arrays for $\Gamma$ are $\{15,8,1;1,8,15\}$, $\{25,12,1;1,12,25\}$ and
$\{35,16,1;1,16,35\}$. Since graphs with these intersection arrays are Taylor graphs, it is impossible by Theorem \ref{general-taylor}. This completes the proof. \epf

\begin{remark}\label{geo-rmk}
It seems to be difficult to classify distance-regular graphs with $\theta_D=-3$ satisfying $k=3(a_1+1)$ 
and $c_2=1$ (see \cite{gdrg-3, yamazaki}). It follows by \cite{gdrg-3}, \cite{non-exist} 
and \cite{yamazaki} that any distance-regular graph in Theorem \ref{geo-thm} (9) 
(i.e., a distance-regular graph with $\theta_D=-3$ satisfying $k=3(a_1+1)$, $c_2=1$ and 
$a_1,\frac{1}{3}c_D\geq 2$) is one of the following:\\
(1) A distance-regular graph satisfying
$D=\HH+2\geq 4$ and $$(c_i,a_i,b_i)=\left \{ \begin{array}{ll}
(1,\alpha,2\alpha+2)& \mbox{ for }1\leq i\leq \HH \\
(2,2\alpha+\beta-1,\alpha-\beta+2) & \mbox{ for }i= \HH +1,\\
(3\beta,3\alpha-3\beta+3,0) & \mbox{ for }i= \HH +2
\end{array}
\right.\mbox{where~} \alpha, \beta \geq 2;$$

(2) A distance-regular graph satisfying
$D=\HH + 2\geq 3$ and $$(c_i,a_i,b_i)=\left\{ \begin{array}{ll}
(1,\alpha,2\alpha+2)& \mbox{ for }1\leq i\leq \HH \\
(1,\alpha+2\beta-2,2\alpha-2\beta+4) & \mbox{ for }i=\HH +1,\\
(3\beta,3\alpha-3\beta+3,0) & \mbox{ for }i= \HH +2
\end{array}
\right.\mbox{where~} \alpha, \beta \geq 2;$$

(3) A halved graph of a distance-biregular graph with vertices of valency $3$ and
$$(c_i,a_i,b_i)=\left\{ \begin{array}{ll}
(1,\alpha,2\alpha+2)& \mbox{ for }1\leq i\leq \HH\\
(1,\alpha+2,2\alpha)& \mbox{ for }i= \HH +1\\
(4,2\alpha-1,\alpha)& \mbox{ for } \HH +2\leq i\leq D-2,\\
(4,2\alpha+\beta-3,\alpha-\beta+2) & \mbox{ for }i=D-1\\
(3\beta,3\alpha-3\beta+3,0) & \mbox{ for }i=D
\end{array}\right.\mbox{where~} \alpha \geq 2 \mbox{~and~} \beta \in \{2,3\}.$$

It is unknown whether there exists a constant $Y$ such that any distance-regular graph 
with $k=3(a_1+1)$ and $c_2=1$ has diameter $D \leq Y$. 
\end{remark}

\begin{quest}
Is it true that any geometric distance-regular graph  
with smallest eigenvalue $-t$ has a $t$-claw? 
\end{quest}

\section{Non-geometric distance-regular graphs without $4$-claws}\label{sect-nongeom}
In this section, we classify non-geometric distance-regular graphs with diameter $D\geq 3$ 
and valency $k>3$ but without $4$-claws (see Theorem \ref{non-geoandk}). 

\begin{theorem}\label{non-geoandk}
Let $\Gamma$ be a non-geometric distance-regular graph with diameter $D\geq 3$, 
valency $k>3$ but without $4$-claws.
Then $\Gamma$ is one of the following graphs:
\begin{itemize}
\item[(1)] the icosahedron with $\iota(\Gamma)=\{5,2,1;1,2,5\}$;
\item[(2)] the Klein graph with $\iota(\Gamma)=\{7,4,1;1,2,7\}$; 
\item[(3)] the halved $e$-cube with 
$\iota(\Gamma)=\{{e\choose 2},{e-2\choose 2},1;1,6,15\}$ with $e\in \{6,7\}$;
\item[(4)] the Taylor graph with $\iota(\Gamma)=\{k,\frac{k-1}{2},1;1,\frac{k-1}{2},k\}$ with $k\in \{13,17\}$; 
\item[(5)] the Gosset graph with $\iota(\Gamma)=\{27,10,1;1,10,27\}$;
\end{itemize}
or $\Gamma$ has one of the following intersection arrays:
\begin{itemize}
\item[(6)] $\iota(\Gamma)=\{44,24,1;1,12,44\}$;
\item[(7)] $\iota(\Gamma)=\{64,34,1;1,17,64\}$; 
\item[(8)] $\iota(\Gamma)=\{104,54,1;1,27,104\}$.
\end{itemize}
\end{theorem}

To prove Theorem \ref{non-geoandk}, we first show in Subsection \ref{sect-valencybounds}
that the valency $k$ of $\Gamma$ is bounded above by some constant 
(see Proposition \ref{prop-bounds}).
In Subsection \ref{sect-comp}, with the aid of computer we find all 
feasible intersection arrays for distance-regular graphs 
with $D\geq 3$ and $3$-claws but without $4$-claws whose valencies satisfy 
the bounds in Subsection \ref{sect-valencybounds}.
Finally, in Subsection \ref{sect-nonexist}, we show that 
for some of the intersection arrays we found in Subsection \ref{sect-comp}, 
there are no corresponding distance-regular graphs without $4$-claws.
The proof of Theorem \ref{non-geoandk} is given in Subsection \ref{sect-nonexist}.

We first need the following lemma.

\begin{lemma}\label{non-geom-valency-bound}
If $\Gamma$ is a non-geometric distance-regular graph with diameter $D\geq 3$, 
valency $k>3$, and $3$-claws but without $4$-claws,
then $c_2\geq 2$ and $k\leq \frac{8}{3}(a_1+1)$ hold. 
Moreover, if $k\leq 2(a_1+1)$, then $\Gamma$ is one 
of the following graphs:
\begin{itemize}
\item[(1)] the halved $e$-cube with 
$\iota(\Gamma)=\{{e\choose 2},{e-2\choose 2},1;1,6,15\}$ with $e\in \{6,7\}$;
\item[(2)] the Taylor graph with $\iota(\Gamma)=\{k,\frac{k-1}{2},1;1,\frac{k-1}{2},k\}$ with $k\in \{13,17\}$;
\item[(3)] the Gosset graph with $\iota(\Gamma)=\{27,10,1;1,10,27\}$.
\end{itemize}
\end{lemma}
\pf If $c_2=1$, then $k=3(a_1+1)$ and 
$\Gamma$ is geometric by \cite[Theorem~3.2]{gdrg-3}.
As $\Gamma$ is non-geometric and $k>3$, it follows from Proposition \ref{prop-large-valency} that $k\leq \frac{8}{3}(a_1+1)$ holds. 
Suppose that $k\leq 2(a_1+1)$ holds. 
By Proposition \ref{prop-a1-half-valency}, the graph $\Gamma$ is the halved $7$-cube or a Taylor graph. 
In the latter case, the result follows from Theorem \ref{general-taylor}.\epf

For the rest of this section, let $\Gamma$ be a non-geometric distance-regular graph with 
diameter $D\geq 3$, valency $k>3$, and $3$-claws but without $4$-claws.
By Lemma \ref{non-geom-valency-bound}, it follows that $k\leq \frac{8}{3}(a_1+1)$, 
and, if $k\leq 2(a_1+1)$, then $k\leq 27$. 
Thus, in what follows, we assume that the valency $k$ satisfies 
\begin{equation}\label{eq-gammadelta-bound}
k=\delta (a_1+1)\text{~where~}2<\delta\leq \frac{8}{3}.
\end{equation}

Further, let $\theta_0=k$ $>$ $\theta_1$ $>$ $\ldots$ $>$ $\theta_D$ be 
the $D+1$ distinct eigenvalues of $\Gamma$ with multiplicities $m_0=1$, 
$m_1$, $\ldots$, $m_D$, respectively.
Lemma \ref{vka} (3) shows that $D\leq 4$.
By Lemma \ref{terw-ev-lemma} (1) (with $t=3$ and $2<\delta\leq \frac{8}{3}$), $\theta_D>-6$ holds.

\subsection{Valency bounds}\label{sect-valencybounds}

It is known that for given integer $t\geq 2$, 
there are only finitely many non-geometric distance-regular graphs 
with both valency and diameter at least $3$ and with smallest eigenvalue at least $-t$, 
see \cite[Theorem~9.10]{drgsurvey}. 
However, in order to perform a computer search of feasible intersection arrays, 
we need to improve the valency bounds.

In this subsection we obtain a better bound for valency $k$
of a non-geometric distance-regular graph $\Gamma$ with diameter $D\geq 3$ 
but without $4$-claws.

To obtain a valency bound, we consider the following two cases: $\theta_1\geq a_1$ and $\theta_1<a_1$. 
In the latter case, we also consider several subcases depending 
on ($m_1\geq k$ or $m_1<k$) and ($\theta_D\in \mathbb{Z}$ or $\theta_D\notin \mathbb{Z}$).

The main result of this subsection is the following proposition.

\begin{prop}\label{prop-bounds}
Let $\Gamma$ be a non-geometric distance-regular graph with diameter $D\geq 3$, 
valency $k>3$, and $3$-claws but without $4$-claws. 
Then one of the following holds.
\begin{enumerate}
\item[(1)] $\theta_1\geq a_1$ holds, $D\in \{3,4\}$, and: 
\begin{enumerate}
\item[(i)] if $\theta_1\in \mathbb{Z}$,
then $k\leq 476$ and $m_1\geq k$,
\item[(ii)] if $\theta_1\notin \mathbb{Z}$, 
then $k\leq 952$ and $m_1\geq \frac{k}{2}$.
\end{enumerate}
\item[(2)] $\theta_1<a_1$ holds, $D=3$, and: 
\begin{enumerate}
\item[(i)] if $m_1\geq k$, then $k\leq 790$,
\item[(ii)] if $m_1<k$ and $\theta_1\notin \mathbb{Z}$, 
then $k\leq 36$,
\item[(iii)] if $m_1<k$, $\theta_1\in \mathbb{Z}$ 
and $\theta_3\notin\mathbb{Z}$, then $k\leq 530$,
\item[(iv)] if $m_1<k$, $\theta_1\in \mathbb{Z}$ 
and $\theta_3\in\mathbb{Z}$, then $k\leq 833$.
\end{enumerate}
\end{enumerate}
\end{prop}

To prove Proposition \ref{prop-bounds}, we need several lemmas.

\begin{lemma}\label{lemma-big-and-int-th1}
If $\theta_1\in \mathbb{Z}$ and $\theta_1\geq a_1$, then $m_1\geq k$. 
\end{lemma}
\pf If $m_1<k$ and $\theta_1$ is integral, then it follows by \cite[Theorem 4.4.4]{bcn}
that $\theta_1+1$ divides $b_1$, i.e., $b_1=\ell (\theta_1+1)$  
holds for some integer $\ell \geq 1$.

By Eq. (\ref{eq-gammadelta-bound}), $b_1\leq \frac{5}{3}(a_1+1)$. 
As $\ell (a_1+1)\leq \ell(\theta_1+1)=b_1 \leq \frac{5}{3}(a_1+1)$,
we find $\ell=1$ and thus $\theta_1=b_1-1$. 
As $\Gamma$ is non-geometric and it satisfies $D\geq 3$, $\theta_1\geq a_1$ with $m_1 < k$,
and $c_2\geq 2$ (see Lemma \ref{non-geom-valency-bound}), 
it follows by \cite[Theorem 4.4.11]{bcn} that
there are no such graphs. Thus $m_1\geq k$, 
and the lemma is proved.
\epf

\begin{lemma} \label{theta1>a1}
If $\theta_1\geq a_1$ and $m_1\geq \frac{k}{N}$ both hold for some constant $N \geq 1$, then  
$$k < \frac{4288}{9}N.$$
\end{lemma}
\pf Suppose that $\theta_1\geq a_1$ and $m_1\geq \frac{k}{N}$ hold for some constant $N\geq 1$. 
It follows by Eq. (\ref{mi}), $m_1\geq \frac{k}{N}$ and $\theta_1\geq a_1$ that $$N v \geq k\left( \sum_{i=0}^{D}k_i u_i^2(\theta_1) \right) > k(1+ku_1^2(\theta_1))= k\left( 1+\frac{\theta_1^2}{k} \right )\geq k+a_1^2$$
and hence
\begin{equation}\label{lem-m1}
v>\frac{k+a_1^2}{N}.
\end{equation}

We also have that:
\begin{itemize}
\item $b_1=k-(a_1+1)=(\delta-1)(a_1+1)$ by Eq. (\ref{eq-gammadelta-bound}),
\item $b_3\leq b_2\leq (a_1+1)$ by Lemma \ref{vka} (2),
\item $c_2\geq \frac{a_1-6}{6}$ holds by Lemma \ref{surveyc2} as $\Gamma$ is not geometric,
\item $c_3\geq b_3+1-b_1+2(a_1+2)\geq 1-b_1+2(a_1+2)>-k+3(a_1+1)=(3-\delta)(a_1+1)$ by 
Eq. (\ref{eq-gammadelta-bound}) and Lemmas \ref{nonterwkbd} and \ref{vka} (2), 
\item if $D=4$, then $c_4\geq 1-b_1+3(a_1+2)>-k+4(a_1+1)=(4-\delta)(a_1+1)$ 
by Eq. (\ref{eq-gammadelta-bound}) and Lemmas \ref{nonterwkbd} and \ref{vka} (2).
\end{itemize}

Suppose $a_1\geq 7$ so that $a_1-6>0$. Then we find by Eq. (\ref{ki}) that
\begin{equation*}
k_2\leq 6(\delta -1)\left( \frac{a_1+1}{a_1-6}\right) k,~~k_3< \frac{6(\delta -1)}{3-\delta} \left( \frac{a_1+1}{a_1-6}\right) k,~\text{and,~if~}D=4,~~
k_4< \frac{6(\delta -1)}{(3-\delta)(4-\delta)}\left( \frac{a_1+1}{a_1-6}\right) k,
\end{equation*}
and hence
\begin{equation*}\label{nongeo-bigm-vbd}
v<1 + k + 6k\left( \frac{a_1+1}{a_1-6}\right)
\left(\frac{(\delta^2-8\delta+17)(\delta-1)}{(3-\delta)(4-\delta)}\right).
\end{equation*}

Note that the function $\frac{(\delta^2-8\delta+17)(\delta-1)}{(3-\delta)(4-\delta)}$ 
is strictly increasing on $2<\delta\leq \frac{8}{3}$, therefore 
\begin{equation*}
v<1 + k + k\frac{125(a_1+1)}{2(a_1-6)}.
\end{equation*}

Suppose $a_1\geq 146$. Then $\frac{(a_1+1)}{(a_1-6)}\leq \frac{21}{20}$ and 
thus $v<67k$. It follows by Eq. (\ref{lem-m1}) that $a_1^2<(67N-1) k \leq \frac{8(67N-1)}{3} (a_1+1)$
and thus
\[a_1< \frac{4(67N-1)+2 \sqrt{2(67N-1)(134N+1)}}{3} \mbox{~and~hence~}k<\max\left\{ \frac{8}{3}(146+1), 
\frac{4288}{9}N \right\}= \frac{4288}{9}N,\]
which completes the proof. \epf

The following technical lemma will be used later (see Lemma \ref{smallthetabigm1}).

\begin{lemma}\label{lem-m}
If $\theta_1<a_1$ and $m_1\geq \frac{k}{N}$ for some constant $N\geq 1$, then $D=3$ 
and 
$$\frac{k b_1}{-N(\theta_3+1)} < ( - \theta_3)\left( v-1-\frac{N+1}{N}k\right) +\frac{1}{N}k.$$
\end{lemma}
\pf Suppose that $\theta_1<a_1$ and $m_1\geq \frac{1}{N}k$ hold for some constant $N\geq 1$. 
By Lemma \ref{d=3}, $D=3$ holds. Using Eq. (\ref{trace}) we have
\begin{equation}\label{eq5.5}
v=tr(A^0)+tr(A)=k+1+(\theta_1+1)m_1+(\theta_2+1) m_2+
(\theta_3+1) m_3.
\end{equation}
Note that $m_2 \geq k$ by \cite[Theorem 4.4.4]{bcn} and thus $m_3=v-1-m_1-m_2 \leq v-1-\frac{N+1}{N}k$. 

By Lemma \ref{vka} (4), $\theta_2\geq -1$ holds and hence 
$(\theta_2+1)m_2\geq 0$. It follows by Eq. (\ref{eq5.5}) that
$$(\theta_1+1)m_1=v-\sum_{i\in \{0,2,3\}} (\theta_i+1) m_i\leq v-k-1 - (\theta_3+1)\left( v-1-\frac{N+1}{N}k \right)=-\theta_3 \left( v-1-\frac{N+1}{N}k \right)+\frac{1}{N}k.$$
Now the result follows as 
$m_1 (\theta_1+1) \geq \frac{1}{N}k (\theta_1+1) > \frac{k b_1}{-N(\theta_3+1)}$ 
holds by Theorem \ref{kypeq}.
\epf

\begin{lemma}\label{smallthetabigm1}
If $\theta_1<a_1$ and $m_1 \geq k$ both hold, then $k\leq 790$.
\end{lemma}
\pf By Eq. (\ref{eq-gammadelta-bound}) and Lemmas \ref{d=3}, \ref{vka} (2) and \ref{surveyc2}, 
we have that
\begin{equation} \label{abcbound}
D=3,~b_1=(\delta-1)(a_1+1)<c_3,~b_2\leq a_1+1,\text{~and}~c_2\geq \frac{a_1-6}{6}.
\end{equation}

\begin{claim}\label{bigm1smalltheta1bound}
The valency $k$ satisfies
\[k < \max \left\{\frac{\delta(\delta +3)(18 \delta^2 +5 \delta -3)}{(\delta -1) (2 \delta -3)^2}, 
(42 \delta -34)\delta, 50 \delta \right \}.\]
\end{claim}
\noindent{\em Proof of Claim \ref{bigm1smalltheta1bound}: }
If $a_1< \max\{42 \delta-35, 49\}$, then $k < \max \left\{(42 \delta -34)\delta , 50 \delta \right \}$. 
Now we assume $a_1\geq  \max\{42 \delta-35, 49\}$. We first show $v<(6\delta +3)k$. 
By Eq. (\ref{abcbound}), we find
\[k_2\leq \frac{k (\delta-1)(a_1+1)}{\frac{a_1-6}{6}}<(6\delta -5)k-1,~~~ k_3\leq \frac{b_1(a_1+1)k}{\left( \frac{a_1-6}{6}\right) c_3}<
\frac{6(a_1+1)}{a_1-6} k <7k \]
and thus $v=1+k+k_2+k_3<(6\delta +3)k$. 

As $\theta_3>-\frac{3\delta}{2\delta -3}$ holds by
Lemma \ref{terw-ev-lemma}, we find $\theta_1+1 > \frac{(2\delta -3) (k-1-a_1)}{\delta+3}$ by Theorem \ref{kypeq}. It follows by Lemma \ref{lem-m}
that
\[\frac{k (2\delta -3)(k-1-a_1)}{\delta +3}<\left( \frac{3\delta}{2\delta-3}\right) (6\delta+1) k+k \]
and thus $k<(a_1+1)+\frac{(18 \delta^2+5\delta-3)(\delta+3)}{(2\delta -3)^2}$. 
As $\delta (a_1+1)=k$, we find $a_1+1<
\frac{(18 \delta^2+5\delta-3)(\delta+3)}{(\delta -1)(2\delta -3)^2}$ and hence 
$k=\delta (a_1+1)<\frac{\delta(18 \delta^2+5\delta-3)(\delta+3)}{(\delta-1)(2\delta -3)^2}$. 
Now the claim follows.\epf

Evaluating the valency bound from Claim \ref{bigm1smalltheta1bound} for $2<\delta\leq \frac{8}{3}$ 
shows that $k\leq 790$ holds, and this completes the proof. \epf

\begin{remark}
Lemma \ref{m1ineq} together with Lemma \ref{terw-ev-lemma} show 
that the multiplicity $m_1$ can be bounded below by $m_1\geq \frac{1}{N}k$ 
with an appropriate constant $N>1$. Therefore the valency $k$ can 
be bounded above in the same manner as in the proof of Lemma \ref{smallthetabigm1}.
However, in the following several lemmas, we will obtain more accurate bounds for $k$ 
for the case when $\theta_1<a_1$ and $m_1<k$. 
\end{remark}

\begin{lemma}\label{theta1-notinteger}
Suppose that $\theta_1<a_1$, $\theta_1\not\in \mathbb{Z}$, and $m_1<k$.
Then $k\leq 36$.
\end{lemma}
\proof 
It follows from Lemma \ref{d=3} and \cite[Theorem~4.4.4]{bcn} that $D=3$ and $m_1=m_3<k\leq m_2$, 
and thus, by Lemma \ref{vka} (4), $\theta_2$ is integral satisfying $\theta_2\geq -1$.

We first show that if $\theta_2=-1$, then $k\leq 36$. 
If $c_2=1$, then $k\leq 36$ holds by \cite[Theorem~4.2]{bangkoolen-3}. 
Let $c_2\geq 2$. It follows by Eq. (\ref{eq-gammadelta-bound}) and \cite[Theorem~4.2~(iv)]{bangkoolen-3} 
that
\begin{equation}\label{-6kbd}
2(\alpha+2)\leq 2(a_1+1){<}k\leq \alpha+11+\frac{50}{\alpha+2}
\end{equation}
where ${\alpha}:=a_1-c_2$. 
{It follows from Eq. (\ref{-6kbd}) that $0\leq \alpha \leq 10$, and thus $k\leq 36$ holds, 
if $\theta_2=-1$.}

Now we assume that $\theta_2\geq 0$. We will show that 
$k\leq 33$. 
If follows from Eq. (\ref{trace}), $m_2\geq k > m_1=m_3$ and $\theta_2\geq 0$ that 
\begin{equation*}
0=\frac{1}{k}tr(A)=1+\frac{m_1}{k}(\theta_1+\theta_3)+\frac{m_2}{k}\theta_2\geq 
1+\frac{m_1}{k}(\theta_1+\theta_3)+\theta_2,
\end{equation*}
and hence $\theta_1+\theta_3\leq \frac{k}{m_1}(-1-\theta_2)<-1-\theta_2$ and 
$\theta_1<-\theta_3-1-\theta_2\leq 5$. Thus, we have that $\theta_1<5$ and $|\theta_3|\leq 6$.

If $k\geq 31>\frac{5^2+6^2}{2}> \frac{\theta_1^2+\theta_3^2}{2}$, then 
it follows by {$tr(A^0)=v$,} $tr(A^2){=vk}$ and $m_1\geq \frac{k}{2}$ (see Lemma 2.6 (2)) 
that 
\begin{equation*}
k(k-1)=m_1(2k-\theta_1^2-\theta_3^2)+m_2(k-\theta_2^2)\geq m_1(2k-\theta_1^2-\theta_3^2)+k(k-4)> \frac{k}{2}(2k-\theta_1^2-\theta_3^2)+k(k-4).
\end{equation*}

Hence $2k-6\leq \theta_1^2+\theta_3^2<25+36$.  
This shows that $k\leq 33$ if $\theta_2\geq 0$. 
This completes the proof. \epf

\begin{lemma}\label{vestimation}
If $\theta_1\in \mathbb{Z}$ and $m_1 <k$, 
then $D=3$ and $b_1=\ell (\theta_1+1)$ hold for some $\ell \in \{2,3,4\}$. 
Moreover, 
$$ v <1 + \delta(a_1+1) + \frac{\delta(\delta-1)\ell}{(\ell-\delta+1)}(a_1+1) +
\frac{\delta^2 \ell }{(\ell-\delta+1)(-\theta_3)}(a_1+1)<
1+k+ \frac{2(\delta-1) \ell}{(\ell-\delta+1)}k.$$
\end{lemma}
\proof 
Suppose that $\theta_1$ is integral with multiplicity $m_1<k$. 
By Lemma \ref{lemma-big-and-int-th1}, we have $\theta_1<a_1$, 
and, by Lemma \ref{d=3}, we have $D=3$.

By \cite[Theorem 4.4.4]{bcn}, $b_1=\ell(\theta_1+1)$ holds for some integer $\ell \geq 1$, 
and $-1-\frac{b_1}{\theta_1+1}=-1-\ell$ is an eigenvalue of the local graph 
$\langle\Gamma_1(x)\rangle$ for a vertex $x\in V(\Gamma)$.
As $k>2(a_1+1)$ holds by Eq. (\ref{eq-gammadelta-bound}), 
it follows by Lemma \ref{terw-ev-lemma} that $\theta_3>-6$ 
and thus $\ell\in \{1,2,3,4\}$ as $-1-\frac{b_1}{\theta_1+1}\geq \theta_3$ follows 
by the interlacing result (see \cite[Corollary 3.3.2]{bcn}).
Note that if $\ell=1$, then by \cite[Theorem 4.4.11]{bcn} $\Gamma$ is 
the halved $e$-cube with $e \in \{6,7\}$, or the Gosset graph, 
i.e., $\Gamma$ satisfies $k\leq 2(a_1+1)$, a contradiction.

Suppose $\ell\in \{2,3,4\}$. By Eq. (\ref{eq-gammadelta-bound}), $b_1=(\delta-1)(a_1+1)$ holds. 
By Lemma \ref{KYP} (1), we have $b_1=\ell (\theta_1+1) > \ell \theta_1 > \ell (a_1+1-c_2)$ 
and thus
$$ c_2>\frac{\ell-\delta+1}{\ell}(a_1+1) \mbox{~~and~~} k_2=\frac{kb_1}{c_2}<\frac{\delta(\delta-1)\ell}{(\ell-\delta+1)}(a_1+1).$$
(Note here that $\ell -\delta + 1 >0$ holds from $b_1=(\delta-1)(a_1+1)=\ell(\theta_1+1)<\ell (a_1+1)$.) 

By Lemma \ref{KYP} (2) with $\theta_1<a_1$, we have $a_3\leq \theta_1<a_1$, 
and thus $c_3=k-a_3\geq k-\theta_1>k-(a_1+1)=b_1=(\delta-1)(a_1+1)$. 
By Lemma \ref{vka}, we have
\begin{equation}\label{b2-less}
b_2\leq \min\{a_1+1,~a_3+1\}\text{~and~} b_2\leq \frac{k}{-\theta_3}.
\end{equation}

Hence it follows by $c_3>b_1$ and Eq. (\ref{b2-less}) that
$$k_3=\frac{k b_1b_2}{c_2c_3}<\frac{\delta^2 \ell }{(\ell-\delta+1)(-\theta_3)}(a_1+1).$$

Now the result follows as $v=1+k+k_2+k_3$. In particular, using Eq. (\ref{b2-less}) and 
$a_3\leq \theta_1<a_1$,
$$\frac{b_2}{c_3}\leq \frac{a_3+1}{k-a_3}< \frac{a_1}{k-a_1}< 1 \mbox{~~holds~~and~~thus~~}k_3=\frac{k_2 b_2}{c_3}<k_2.$$
This shows $v=1+k+k_2+k_3 < 1+k+2k_2 <1+k+ \frac{2(\delta-1)\ell}{(\ell-\delta+1)} k$. \epf

\begin{lemma}\label{smalltheta1smallm1lastnoninteger}
Suppose that $\theta_1<a_1$, $\theta_1\in \mathbb{Z}$, and $m_1 <k$.
If the smallest eigenvalue of $\Gamma$ is not integral, then $k\leq 530$.
\end{lemma}
\proof
By Lemma \ref{d=3}, we have $D=3$. 
Suppose that $\theta_3\notin \mathbb{Z}$.
Then the eigenvalues $\theta_2$ and $\theta_3$
are conjugate algebraic integers with $m_2=m_3$ as $\theta_1\in \mathbb{Z}$. 
By Eq. (\ref{trace}), we have 
$$tr(A^0)=1+m_1+m_2+m_3=1+m_1+2m_2=v \mbox{~~and~~} tr(A^1)=k+\theta_1m_1+m_2(\theta_2+\theta_3)=0.$$
This shows
\begin{equation}\label{v-lower-bd}
v=1+m_1+ \frac{2(k+\theta_1m_1)}{-(\theta_2+\theta_3)}.
\end{equation}

By Lemma \ref{vestimation}, $b_1=\ell(\theta_1+1)=(\delta-1)(a_1+1)$ holds for some $\ell\in \{2,3,4\}$.
In particular,
\begin{equation}\label{eq-th1-lower-bd}
\theta_1=\frac{(\delta-1)(a_1+1)}{\ell}-1\ge \frac{\delta-1-\epsilon}{\ell } (a_1+1)
\end{equation}
holds for some $\frac{\ell}{a_1+1}\leq \epsilon \leq \delta-1$. 
If $a_1+1< 200$, then $k\leq 530$ as $k\leq \frac{8}{3}(a_1+1)$. 

In the rest of the proof,
we assume $a_1+1 \geq 200$ and $\epsilon=0.02$. 
It follows by Lemma \ref{vka} (4) that $0< -(\theta_2+\theta_3)\le 1-\theta_3$.
We find by Eqs. (\ref{v-lower-bd}) and (\ref{eq-th1-lower-bd}) that
\begin{equation}\label{v-2nd-lower-bd}
v=1+m_1+ \frac{2(k+\theta_1m_1)}{-(\theta_2+\theta_3)}> 1+
\frac{2 m_1(\delta-1.02)}{(1-\theta_3)\ell}(a_1+1).
\end{equation}
By Lemma \ref{vestimation} and Eq. (\ref{v-2nd-lower-bd}),
$$ \frac{2 (\delta-1.02)}{(1-\theta_3)\ell} m_1<
\delta + \frac{\delta(\delta-1) \ell}{(\ell-\delta+1)} + \frac{\frac{\delta^2 \ell}{-
\theta_3}}{(\ell-\delta+1) }$$
follows and thus
\[m_1 < \frac{\ell (1-\theta_3)}{2}\times \frac{\delta^2(\ell-1) + \delta + \frac{\delta^2 \ell}{-\theta_3}}{(\delta-1.02)(\ell-\delta+1)}.\]
We now put $ M(\delta,\theta_3,\ell):=\frac{\ell (1-\theta_3)}{2}\times \frac{A(\delta,\theta_3,\ell)}{B(\delta,\ell)}$, where
\begin{equation}\label{AB}
A(\delta,\theta_3,\ell):=\delta^2(\ell-1) + \delta + \frac{\delta^2 \ell}{-\theta_3} \mbox{~~and~~}B(\delta,\ell):=(\delta-1.02)(\ell-\delta+1).
\end{equation}
Then
\begin{equation}\label{m1-upper-bd}
m_1<M(\delta,\theta_3, \ell).
\end{equation}
To determine an upper bound in Eq. (\ref{m1-upper-bd}), we first show the following claim.

\begin{claim}\label{small-theta1-claim} Let $\tau$ and $\gamma$ be some real numbers.\\
(1) If $\tau \leq \theta_3\leq -2$ holds, then $M(\delta,\theta_3,\ell)\leq M(\delta,\tau,\ell)$.\\
(2) If $2< \delta\leq \gamma\leq \frac{8}{3}$ holds, then 
$M(\delta,\theta_3,\ell)\le M(\gamma,\theta_3,\ell)$.
\end{claim}
\noindent{\em Proof of Claim \ref{small-theta1-claim}: }
(1) For given $\delta$ and $\ell$, function $A(\delta,\theta_3,\ell)(1-\theta_3)$ of $\theta_3$ is decreasing on $\theta_3\leq -2$ as
$$ \frac{\partial}{\partial  \theta_3}\left\{ A(\delta,\theta_3,\ell)(1-\theta_3)\right\}= \frac{\ell \delta^2}{\theta_3^2}-\delta(1 + (\ell-1)\delta)
=\delta \left(\frac{\ell \delta}{\theta_3^2}-1-(\ell-1)\delta\right)\leq
\delta \left( \frac{\ell\delta}{4}-1-(\ell-1)\delta \right)<0.$$
Thus $M(\delta,\theta_3,\ell)$ is also decreasing on $\theta_3\leq -2$, and this shows part (1) as $\tau \leq \theta_3\leq -2$. \\
\noindent (2) As $ M(\delta,\theta_3,\ell)=\frac{\ell (1-\theta_3)}{2}\times \frac{A(\delta,\theta_3,\ell)}{B(\delta,\ell)}$,
it is enough to show that $\frac{\partial}{\partial \delta}\left( \frac{A(\delta,\theta_3,\ell)}{B(\delta,\ell)} \right)>0$ holds for all
$2< \delta \leq \frac{8}{3}$. Note that $2<\delta <3\leq \ell+1$. By Eq. (\ref{AB}), we have
\begin{eqnarray}\label{partial1}
& B(\delta,\ell)^2 \times \frac{\partial}{\partial \delta}\left( \frac{A(\delta,\theta_3,\ell)}{B(\delta,\ell)} \right)
=\frac{\partial A(\delta,\theta_3,\ell) }{\partial \delta} B(\delta,\ell)- A(\delta,\theta_3,\ell)\frac{\partial B(\delta,\ell)}{\partial \delta}=\nonumber\\
&=\left((2\ell-2)\delta + 1 + \frac{2\ell\delta}{-\theta_3}\right)(\delta-1.02)(\ell-\delta+1) -
\left(\delta^2(\ell-1) + \delta + \frac{\ell\delta^2}{-\theta_3}\right)(\ell- 2\delta+2.02).
\end{eqnarray}
If $\ell - 2\delta +2.02<0$ then $\frac{\partial}{\partial \delta}\left( \frac{A(\delta,\theta_3,\ell)}{B(\delta,\ell)} \right)
>0$ holds by Eq. (\ref{partial1}). 
Now we suppose that  $\ell - 2\delta +2.02 \geq 0$ holds. Then for each $\ell \in \{2,3,4\}$,
$2< \delta\le \min \left\{ \frac{\ell+2.02}{2}, \frac{8}{3} \right\}$ holds and thus
\begin{eqnarray*}
\left((2\ell-2)\delta + 1\right)(\delta-1.02)(\ell-\delta+1) &>&
\left(\delta^2(\ell-1) + \delta \right)(\ell- 2\delta+2.02),\\
\frac{2\ell\delta}{-\theta_3}(\delta-1.02)(\ell-\delta+1)&>&
\frac{\ell\delta^2}{-\theta_3}(\ell- 2\delta+2.02).
\end{eqnarray*}
This yields that the right hand side in Eq. (\ref{partial1}) is positive. This proves (2) and hence the claim follows. \epf

To complete the proof of this lemma, we consider the following two cases. 

If $\frac{7}{3}<\delta\leq \frac{8}{3}$, then it follows by Lemma \ref{terw-ev-lemma} (1)
that there exists a clique $C$ satisfying $|C|> \frac{5k}{21}+1$ and thus $\theta_{3}>-\frac{21}{5}$.
 By Eq. (\ref{m1-upper-bd}) and Claim \ref{small-theta1-claim}, we find $m_1\leq 124$
as 
\[
m_1<M(\delta,\theta_3, \ell)\leq \max \left\{ M\left( \frac{8}{3},
-\frac{21}{5},\ell \right) \mid \ell \in \{2,3,4\} \right\}=M\left( \frac{8}{3},-\frac{21}{5},2 \right)<125.\]
It follows by Lemma \ref{m1ineq} (1) that $\frac{5k}{21}+1\leq m_1\leq 124$ and thus $k\leq 516$.

Now we suppose that $2<\delta\leq \frac{7}{3}$. By Lemma \ref{terw-ev-lemma} (1), 
there exists a clique $C$ satisfying $|C|> \frac{k}{6}+1$ and thus $\theta_{3}>-6$.
Now we find $m_1\leq 89$ as
\[
m_1<M(\delta,\theta_3, \ell)\leq \max \left\{ M\left( \frac{7}{3},-6,\ell \right) \mid \ell \in \{2,3,4\} \right\}=M\left( \frac{7}{3},-6,4 \right)<90
\]
holds by Eq. (\ref{m1-upper-bd}) and Claim \ref{small-theta1-claim}. 
Hence it follows by Lemma \ref{m1ineq} (1) that $\frac{k}{6}+1\leq m_1\leq 89$ and thus $k\leq 528$.

This shows that if $a_1+1\geq 200$, then $k\leq 528$. 
Recall that $k\leq 530$, if $a_1+1<200$.
This completes the proof.\epf

\begin{lemma}\label{smalltheta1smallm1lastinteger}
Suppose that $\theta_1<a_1$, $\theta_1\in \mathbb{Z}$ and $m_1 <k$.
If the smallest eigenvalue of $\Gamma$ is integral, then $k \leq 833$.
\end{lemma}
\proof 
By Eq. (\ref{eq-gammadelta-bound}), Lemmas \ref{non-geom-valency-bound} and \ref{vestimation}, 
we have that $k=\delta(a_1+1)$ for some number 
$2<\delta\leq \frac{8}{3}$, $b_1=\ell(\theta_1+1)$ for some $\ell\in \{2,3,4\}$, and $D=3$.

Using Eq. (\ref{trace}), we have $v=tr(A^0)+tr(A)=k+1+(\theta_1+1)m_1+(\theta_2+1)m_2+(\theta_3+1)m_3.$ Note that $m_2\geq k$ by \cite[Theorem 4.4.4]{bcn}
and thus $m_3=v-1-m_1-m_2\leq v-2-k$. As $\theta_2\geq -1$ holds by Lemma \ref{vka} (4),
we obtain that 
$$(\theta_1+1)m_1=v-\sum_{i\in \{0,2,3\}} (\theta_i+1)m_i\leq v-k-1 - (\theta_3+1)(v-k-1)=-\theta_3(v-k-1).$$

By Lemma \ref{vestimation}, $b_1=\ell (\theta_1+1)$ and $k=\frac{b_1 \delta}{\delta -1}$, it follows that
\begin{equation}\label{integerm1bd}
m_1\leq \frac{-\theta_3 (v-k-1)}{\theta_1+1} \leq -2\theta_3\frac{\ell^2 \delta}{\ell-\delta+1}.
\end{equation}

To complete the proof of the lemma, we consider the following three cases. 

If $2<\delta\leq \frac73$, then it follows by Lemma \ref{terw-ev-lemma} (1) 
that there exists a clique $C$ satisfying $|C|> \frac{k}{6}+1$ and thus $\theta_{3}\geq -5>-6$. 
Hence it follows by Lemma \ref{m1ineq} (1) and Eq. (\ref{integerm1bd})
with $(\delta, \theta_3,\ell)=(\frac73,-5,2)$ that
$\frac{k}{6}+1 < m_1\leq 140$. This shows $k < 834$.

If $\frac73 <\delta\leq \frac52 $, then it follows by Lemma \ref{terw-ev-lemma} (1) that
there exists a clique $C$ satisfying $|C|> \frac{5k}{21}+1$ and thus $\theta_{3}\geq -4>-\frac{21}{5}$. 
Hence it follows by Lemma \ref{m1ineq} (1) and Eq. (\ref{integerm1bd})
with $(\delta, \theta_3,\ell)=(\frac52,-4,2)$ that
$\frac{5k}{21}+1<  m_1\leq 160$. This shows $k< 668$. 

If $\frac52 <\delta\leq \frac83 $, then it follows by Lemma \ref{terw-ev-lemma} (1) that
there exists a clique $C$ satisfying $|C|> \frac{4k}{15}+1$ and thus $\theta_{3}\geq -3>-\frac{15}{4}$. 
Hence it follows by Lemma \ref{m1ineq} (1) and Eq. (\ref{integerm1bd})
with $(\delta, \theta_3,\ell)=(\frac83,-3,2)$ that
$\frac{4k}{15}+1 < m_1 \leq 192$. This shows $k < 717$. 

This completes the proof. \epf

\noindent{\em Proof of Proposition \ref{prop-bounds}: }
By Lemma \ref{vka} (3), $D\in \{3,4\}$ holds.

(1) Suppose that $\theta_1\geq a_1$ holds. 
If $\theta_1\in \mathbb{Z}$, then $m_1\geq k$ holds by Lemma \ref{lemma-big-and-int-th1}, 
and thus $k\leq 476$ by Lemma \ref{theta1>a1}. 
If $\theta_1\notin \mathbb{Z}$, then $m_1\geq k/2$ holds by Lemma \ref{m1ineq} (2), 
and thus $k\leq 952$ by Lemma \ref{theta1>a1}. This shows (1).

(2) Suppose that $\theta_1<a_1$. Then $D=3$ holds by Lemma \ref{d=3}.
The statements (i)--(iv) follow by Lemmas 
\ref{smallthetabigm1}, \ref{theta1-notinteger}, \ref{smalltheta1smallm1lastnoninteger}, \ref{smalltheta1smallm1lastinteger}, respectively. 
This completes the proof.\epf

\subsection{Computational results}\label{sect-comp}

In this subsection, with the aid of computer, we find all 
feasible intersection arrays for distance-regular graphs 
with $D\geq 3$ and $3$-claws but without $4$-claws whose valency satisfies 
the bounds in Proposition \ref{prop-bounds}.

We call an intersection array $\{b_0,\ldots,b_{D-1};1,c_2,\ldots,c_D\}$ 
{\it feasible}, if it satisfies the following conditions:
\begin{itemize}
\item[$(F1)$] $k_i\in \mathbb{N}~~(i=2,\ldots,D)$;
\item[$(F2)$] $ka_1\equiv \ldots \equiv k_Da_D\equiv 0~({\rm mod}~2)$ (the handshake lemma);
\item[$(F3)$] $vk_ia_i\equiv 0~({\rm mod}~6)$ (see \cite[Lemma 4.3.1]{bcn}),
where $v=|V(\Gamma)|=\sum_{i=0}^Dk_i$;
\item[$(F4)$] all multiplicities $m_i$ calculated by Eq. (\ref{mi}) are positive integers;
\item[$(F5)$] all intersection numbers $p_{ij}^k$ are non-negative integers.
\end{itemize}

Suppose that $\Gamma$ is a non-geometric distance-regular graph 
with diameter $D\geq 3$, valency $k>3$, and $3$-claws but without $4$-claws.

For each of the following cases given by Proposition \ref{prop-bounds}:
\begin{enumerate}
\item[$(C1)$] $\theta_1\geq a_1$, $\theta_1\in \mathbb{Z}$, $m_1\geq k$, $D\in \{3,4\}$, and $k\leq 476$;
\item[$(C2)$] $\theta_1\geq a_1$, $\theta_1\notin \mathbb{Z}$, $m_1\geq \frac{k}{2}$, $D\in \{3,4\}$ 
and $k\leq 952$;
\item[$(C3)$] $\theta_1<a_1$, $m_1\geq k$, $D=3$, and $k\leq 790$;
\item[$(C4)$] $\theta_1<a_1$, $m_1<k$, $\theta_1\notin \mathbb{Z}$, $D=3$, and $k\leq 36$;
\item[$(C5)$] $\theta_1<a_1$, $m_1<k$, $\theta_1\in \mathbb{Z}$, $\theta_3\notin\mathbb{Z}$, $D=3$, 
and $k\leq 530$;
\item[$(C6)$] $\theta_1<a_1$, $m_1<k$, $\theta_1\in \mathbb{Z}$, $\theta_3\in\mathbb{Z}$, $D=3$, 
and $k\leq 833$,
\end{enumerate}
we check by computer all feasible intersection arrays additionally satisfying:
\begin{itemize}
\item $2(a_1+1)<k\leq \frac{8}{3}(a_1+1)$ (see Lemma \ref{non-geom-valency-bound}),
\item $\theta_D>-6$ (by Lemma \ref{terw-ev-lemma} with $t=3$ and $\delta=2$),
\item $c_2 \geq \max \left\{2, \frac{a_1-6}{6}\right\}$ (see Lemma \ref{vka} (1) and Lemma \ref{surveyc2}), and
\item $b_2\leq a_3+1$ (by Lemma \ref{vka} (2)).
\end{itemize}

For a given intersection array, we successively check the properties $(Fi)$, 
$i=1,2,\ldots,5$. If it passes $(F1)$--$(F3)$, we then check $(F4)$ by calculating the eigenvalues and 
their multiplicities. For both cases, $D=3$ and $D=4$, we apply a root-finding algorithm known 
as Brent's method \cite{brent} for the characteristic polynomial of the matrix $L_1$. 
For $D=3$, if we find one eigenvalue then we apply the observation in \cite[p. 130]{bcn} 
in order to calculate two remaining non-trivial eigenvalues. For $D=4$, we look for two roots 
in the intervals $(-6,-2)$ and $[a_1,b_1-1)$, respectively (as $\theta_1\leq b_1-1$ 
by \cite[Proposition 4.4.9(i)]{bcn} and $\theta_1\geq a_1$ by the assumption). 
If two such roots are found, then we proceed as in the case $D=3$. 
Finally, if all the multiplicities $m_i$, $i=1,\ldots,D$ are integral and positive, 
then we check whether $(F5)$ holds.

We can now summarize the computational results.

\begin{prop}\label{prop-feasiblearrays}
The following holds.
\begin{enumerate}
\item[(1)] The only feasible intersection arrays satisfying $(C1)$ are:
\[\{10,6,1;1,3,10\},~\{27,16,4;1,2,24\},~\{45,28,4;1,2,42\},
\text{~and~}\{65,40,16;1,4,50\}.\]
\item[(2)] The only feasible intersection array satisfying $(C2)$ is:
\[\{7,4,1;1,2,7\}.\]
\item[(3)] The only feasible intersection arrays satisfying $(C3)$ are:
\[\{13,8,1;1,4,13\},~\{16,10,1;1,5,16\},\text{~and~}\{39,24,1;1,4,39\}.\]
\item[(4)] There are no feasible intersection arrays satisfying $(C4)$.
\item[(5)] There are no feasible intersection arrays satisfying $(C5)$.
\item[(6)] The only feasible intersection arrays satisfying $(C6)$ are:
\[\{15,8,1;1,8,15\},~\{44,24,1;1,12,44\},~\{64,34,1;1,17,64\},~\{104,54,1;1,27,104\},\]
\[\{125,78,1;1,26,125\},~\{155,96,1;1,32,155\},~\{275,168,1;1,56,275\},\]
\[\{575,348,1;1,116,575\},\text{~and~}\{715,432,25;1,80,675\}.\]
\end{enumerate}
\end{prop}

\subsection{Non-existence of some distance-regular graphs without $4$-claws}\label{sect-nonexist}

In this subsection, we show that 
for some of the intersection arrays found in Subsection \ref{sect-comp} 
(see Proposition \ref{prop-feasiblearrays}) 
there are no corresponding distance-regular graphs without $4$-claws.

\begin{lemma}\label{lemma-nonexist-1}
There are no $4$-claw-free distance-regular graphs with intersection arrays satisfying $(C1)$:
\[\{10,6,1;1,3,10\},~\{27,16,4;1,2,24\},~\{45,28,4;1,2,42\},
\text{~and~}\{65,40,16;1,4,50\}.\]
\end{lemma}
\proof
The intersection array $\{10,6,1;1,3,10\}$ is ruled out by \cite[Proposition 1.10.5]{bcn}, 
and the last three arrays are ruled out by Lemma \ref{vka} (1) and Lemma \ref{nonterwkbd}. 
Indeed, they satisfy $k>2(a_1+1)$ and $k<3a_1-3c_2+7$, and
hence any distance-regular graph with one of these three intersection arrays 
must be a Terwilliger graph by Lemma \ref{nonterwkbd}.
However, this is impossible by Lemma \ref{vka} (1).\epf

\begin{lemma}\label{lemma-nonexist-2}
There are no $4$-claw-free distance-regular graphs with the following intersection arrays satisfying $(C6)$:
\[~\{15,8,1;1,8,15\},
~\{125,78,1;1,26,125\},~\{155,96,1;1,32,155\},\]
\[\{275,168,1;1,56,275\},~\{575,348,1;1,116,575\},\text{~and~}\{715,432,25;1,80,675\}.\]
\end{lemma}
\proof The intersection arrays $\{125,78,1;1,26,125\}$, $\{155,96,1;1,32,155\}$, 
$\{275,168,1;1,56,275\}$, $\{575,348,1;1,116,575\}$ and $\{715,432,25;1,80,675\}$ 
satisfy $\theta_3=-5$ and $\frac{5}{2}(a_1+1)<k\leq \frac{8}{3}(a_1+1)$. 
This is impossible as $\theta_3\geq -3$ follows from Lemma \ref{terw-ev-lemma} (1) 
(see the proof of Lemma \ref{smalltheta1smallm1lastinteger}).
 
A distance-regular graph $\Gamma$ with intersection array $\{15,8,1;1,8,15\}$ has 
eigenvalues $\theta_1=3$ and $\theta_3=-5$ with multiplicities $m_1=10$ and $m_3=6$ by Eq. (\ref{mi}). 
It follows by \cite[Theorem 4.4.4]{bcn} that each local graph of $\Gamma$ has eigenvalues $-1-\frac{8}{-5+1}=1$ with multiplicity $15-6 = 9$ 
and $-1-\frac{8}{3+1} = -3$ with multiplicity $15-10 = 5$. Hence each local graph is a strongly 
regular graph with parameters $(15,6,1,3)$, i.e. the unique generalized quadrangle $GQ(2,2)$. 
Since $GQ(2,2)$ is the complement of the triangular graph $T(6)$, which has a 4-clique, 
$GQ(2,2)$ contains a 4-coclique and thus $\Gamma$ has a 4-claw. 
This is impossible.\epf


\begin{remark}
We note that there exists a distance-regular graph with intersection array $\{44,24,1;1,12,44\}$. 
It was constructed in \cite{Klin-Pech}.
One can check that this graph contains $4$-claws, however, 
it is an open question whether a distance-regular graph with this intersection array 
must have a $4$-claw or not.

In the rest of this subsection, we study this problem for the intersection arrays 
satisfying $(C3)$.
We note that there exists a distance-regular graph with intersection array $\{13,8,1;1,4,13\}$, 
given by the Mathon construction, see \cite[Proposition~12.5.3]{bcn}. 
There also exists a distance-regular graph $\Gamma$ with intersection array 
$\{16,10,1;1,5,16\}$, which is locally folded 5-cube, see \cite[p. 386, Remark~(iii)]{bcn}. 
\end{remark}

\begin{lemma}\label{lemma-nonexist-3}
Any distance-regular graph with intersection array 
$\{13,8,1;1,4,13\}$ or $\{16,10,1;1,5,16\}$ contains a $4$-claw.
\end{lemma}
\pf Suppose that $\Gamma$ is a $4$-claw-free distance-regular graph 
with intersection array $\{13,8,1;1,4,13\}$ or $\{16,10,1;1,5,16\}$.
Note that $\theta_3=-\sqrt{k}$ for both possibilities. 

The Delsarte bound (see Eq. (\ref{dc-bd})) for $\Gamma$ shows 
that a clique in $\Gamma$ contains at most $\lfloor1+\sqrt{k}\rfloor$ 
vertices.

Pick a vertex $x\in V(\Gamma)$, and let $y_1,y_2$ be two 
non-adjacent vertices of $\Gamma_1(x)$. As $\Gamma$ is $4$-claw-free, 
the set 
$\langle\{x\}\cup\Gamma_1(x)\setminus \big(\{y_1,y_2\}\cup \Gamma_1(y_1)\cup \Gamma_1(y_2)\big)\rangle$
induces a clique of size $1+k-(2(a_1+1-\nu)+\nu)=k-2a_1-1+\nu$, 
where $\nu=|\Gamma_1(x,y_1,y_2)|$. For $(k,a_1)\in \{(13,4),(16,5)\}$, 
we have $k-2a_1-1+\nu\leq \lfloor1+\sqrt{k}\rfloor$. 
Thus $\nu=0$.


Now pick a vertex $y_3$ from 
$Y:=\Gamma_1(x)\setminus \big(\{y_1,y_2\}\cup \Gamma_1(y_1)\cup \Gamma_1(y_2)\big)$.
As $|\Gamma_1(x,y_i,y_j)|=0$ holds for any $1\leq i<j\leq 3$, 
we see that $k=|\Gamma_1(x)|\geq \sum_{i=1}^3|\{y_i\}\cup \Gamma_1(x,y_i)|=3(a_1+1)$. 
This is impossible, as $k\leq \frac{8}{3}(a_1+1)$, and the lemma follows.\epf

\begin{prop}\label{39}
There are no distance-regular graphs with intersection array $\{39,24,1;1,4,39\}$. 
\end{prop}
\pf Let $\Gamma$ be a distance-regular graph with intersection array $\{39,24,1;1,4,39\}$. 
Then $\theta_0=39$, $\theta_1=13$, $\theta_2=-1$, $\theta_3=-3$, and 
$m_0=1$, $m_1=45$, $m_2=39$, $m_3=195$.

For a vertex $\infty\in V(\Gamma)$, if the local graph $\Delta:=\langle \Gamma_1(\infty) \rangle$ 
contains an $s$-coclique then 
\[
c_2-1\geq \frac{s(a_1+1)-k}{{s \choose 2}}
\]
holds by \cite[Lemma~2]{shilla}. This yields that $\Gamma$ is $4$-claw-free, 
as $4\cdot 15-39>6\cdot 3$.

Let $\{x_1,x_2,x_3\}$ induce a $3$-coclique in $\Delta$. As $\Gamma$ is $4$-claw-free, 
we see that 
\[
39=3\cdot (14+1)-\sum_{1\le i<j\le 3} |\Delta_1(x_i,x_j)|+|\Delta_1(x_1,x_2,x_3)|,
\]
so that $\sum_{1\le i<j\le 3} |\Delta_1(x_i,x_j)|\geq 6$.
This implies that $x_i$ lies in a maximal clique of $\Delta$, say $C_{x_i}$, 
with $\gamma_i \geq 9$ vertices, $i\in\{1,2,3\}$.
The Delsarte bound (see Eq. (\ref{dc-bd})) for $\Gamma$ shows that $\gamma_i\leq 13$.

The number of vertices of $\Delta\setminus \{C_{x_1},C_{x_2},C_{x_3}\}$ is equal to 
\[
\sum_{1\le i<j\le 3} |\Delta_1(x_i,x_j)|-2|\Delta_1(x_1,x_2,x_3)|\leq 6,
\]
so that $\gamma_1+\gamma_2+\gamma_3\geq 33$.

Let ${\cal C}$ be the set of maximal cliques in $\Delta$ with at least $9$ vertices. 
Then $C_{x_1},C_{x_2},C_{x_3}\in {\cal C}$, and for every vertex $x\in \Delta$, 
there exists a clique $C_x\in {\cal C}$, $|C_x|\leq 13$.
As $c_2=4$, we have $|C\cap C'|\leq 3$ for all distinct $C,C'\in {\cal C}$.
Therefore ${\cal C}$ satisfies one of the following:
\begin{itemize}
\item[(a)] ${\cal C}=\{C_1,C_2,C_3\}$, where every $C_i\cup \{\infty\}$ is a Delsarte clique in $\Gamma$;
\item[(b)] ${\cal C}=\{C_1,C_2,C_3,C_4\}$, where $|C_1|=13$, $|C_2|=11$, $|C_3|=|C_4|=9$, 
and $C_3$ and $C_4$ intersect in exactly three vertices;
\item[(c)]${\cal C}=\{C_1,C_2,C_3,C_4\}$, where $|C_1|=12$, $|C_2|=12$, $|C_3|=|C_4|=9$, 
and $C_3$ and $C_4$ intersect in exactly three vertices.
\end{itemize}

\begin{claim}\label{39-b-case}
Case (b) is impossible.
\end{claim}
\noindent{\em Proof of Claim \ref{39-b-case}:}
In this case there exists a vertex $x\in C_4\setminus C_3$ with 6 neighbours in $C_1\cup C_2$ 
($x$ has valency 14, and has no neighbours in $C_3\setminus C_4$), but 
there are at most three of its neighbours in $C_1$ and at most three of its neighbours in $C_2$. 
This implies that $x$ has at least four neighbours in $C_1\cup \{\infty\}$, which is a Delsarte clique 
in $\Gamma$, a contradiction. This shows the claim.\epf

Thus, by Claim \ref{39-b-case} and the connectivity of $\Gamma$, 
for all vertices $x$ in $\Gamma$, the sets ${\cal C}_x$ of maximal cliques 
in $\langle \Gamma_1(x) \rangle$ with at least $9$ vertices 
satisfy either Case $(a)$ or Case $(c)$ from the above.

In Case (a), $\Gamma$ is geometric, a contradiction. 

In Case (c), the number of cliques in $\Gamma$ with exactly $13$ vertices 
is equal to $2v/13=2\cdot 280/13$ and hence is non-integer, a contradiction, 
and the proposition follows. \epf 

We close this section by proving Theorem \ref{non-geoandk}.

\noindent{\em Proof of Theorem \ref{non-geoandk}: }
Let $\Gamma$ be a non-geometric distance-regular graph with diameter 
$D\geq 3$, valency $k>3$ but without $4$-claws.
If $\Gamma$ is $3$-claw-free, then $\Gamma$ is the icosahedron by \cite{no3-claws}.
Suppose that $\Gamma$ contains a $3$-claw. 
If $k\leq 2(a_1+1)$ holds, then, by Lemma \ref{non-geom-valency-bound},
$\Gamma$ satisfies (3), (4), or (5) of Theorem \ref{non-geoandk}.
If $k>2(a_1+1)$, then it follows by Proposition \ref{prop-feasiblearrays}, 
Lemmas \ref{lemma-nonexist-1}, \ref{lemma-nonexist-2}, \ref{lemma-nonexist-3}, 
and Proposition \ref{39} that $\Gamma$ has the intersection array 
$\{7,4,1;1,2,7\}$, i.e., $\Gamma$ is the Klein graph, or one of the intersection arrays 
given in (6)--(8) of Theorem \ref{non-geoandk}. This completes the proof.\epf

\section{Proof of Theorem \ref{theo-main}}

\noindent{\em Proof of Theorem \ref{theo-main}: }
This is straightforward by Proposition \ref{prop-large-valency}, 
Theorems \ref{general-taylor}, Theorem \ref{geo-thm}, and Theorem \ref{non-geoandk}.
\epf


\begin{center}
{\bf Acknowledgements}
\end{center}
Sejeong Bang was supported by the National Research Foundation of Korea (NRF) grant funded by 
the Korea government (MSIP) (No. NRF-2011-0013985).

The research of Alexander Gavrilyuk was funded by Chinese Academy of Sciences President's 
International Fellowship Initiative (Grant No. 2016PE040). His work (e.g., Proposition 
\ref{prop-feasiblearrays}) was also partially supported by the Russian Science Foundation 
(grant 14-11-00061).

Jack Koolen was partially supported by the National Natural Science Foundation of China 
(No. 11471009 and No. 11671376). 
\medskip



\begin{thebibliography}{99}

\bibitem{Babai80}
L. Babai, On the complexity of canonical labeling of strongly regular graphs, 
SIAM J. Computing, 9(1) (1980), 212-216.

\bibitem{BabaiSRG} 
L.~Babai, X.~Chen,~X.~Sun,~S.-H.~Teng,~J.~Wilmes, Faster canonical forms for strongly regular 
graphs. In Proc. 54th Ann. Symp. on Foundations of Computer Science (FOCS'13), 
IEEE Computer Society (2013), 157--166.

\bibitem{Babai}
L.~Babai,~J.~Wilmes, Asymptotic Delsarte cliques in distance-regular graphs, 
J. Alg. Combin. 23(4) (2016), 771--782.

\bibitem{gdrg-3}
S.~Bang, Geometric distance-regular graphs without $4$-claws, Linear Algebra Appl. 
438 (2013), 37--46.

\bibitem{BIConjecture}
S.~Bang,~A.~Dubickas,~J.H.~Koolen,~V.~Moulton, There are only finitely many
distance-regular graphs of fixed valency greater than two, Adv. Math. 269 (2015), 1--55.

\bibitem{DCG1}
S.~Bang, A.~Hiraki,~J.~H.~Koolen, Delsarte clique graphs,
European J. Combin. 28(2) (2007), 501--516.

\bibitem{non-exist}
S.~Bang,~J.H.~Koolen, On geometric distance-regular graphs with diameter three,
European J. Combin. 36 (2014), 331--341.

\bibitem{bangkoolen-3} 
S.~Bang,~J.H.~Koolen, Distance-regular graphs of diameter $3$ having eigenvalue $-1$, 
Linear Algebra Appl. 531 (2017) 38--53.

\bibitem{biggs}
N. Biggs, Algebraic Graph Theory, Second edition, Cambridge University Press, 
Cambridge, 1993.

\bibitem{k=3}
N.L.~Biggs, A.G.~Boshier,~J.~Shawe-Taylor, Cubic distance-regular graphs, 
J. Lond. Math. Soc. (2),  33(3) (1986), 385--394.

\bibitem{no3-claws}
A.~Blokhuis,~A.~E.~Brouwer, Determination of the distance-regular graphs without $3$-claws, 
Discrete Math. 163(1-3) (1997), 225--227.

\bibitem{4clique}
A.~Bondarenko,~A.~Prymak,~D.~Radchenko, Non-existence of $(76,30,8,14)$ strongly regular graph, 
Linear Algebra Appl. 527 (2017) 53--72.

\bibitem{4clique-web}
A.~Bondarenko,~A.~Prymak,~D.~Radchenko, \href{http://prymak.net/SRG-76-30-8-14}{Supplementary files 
for the proof of non-existence of $SRG(76,30,8,14)$}.

\bibitem{Bose}
R.C.~Bose, Strongly regular graphs, partial geometries and partially balanced designs, 
Pacific J. Math. 13 (1963), 389--419.

\bibitem{brent}
R.P.~Brent, Algorithms for Minimization without Derivatives,
Englewood Cliffs, NJ: Prentice-Hall, 1973.

\bibitem{bcn}
A.E.~Brouwer, A.M.~Cohen,~A.~Neumaier, Distance-regular graphs,
Springer-Verlag, Berlin, 1989.

\bibitem{CameronEtAl}
P.J.~Cameron,~J.-M.~Goethals,~J.J.~Seidel,~E.E.~Shult, Line graphs, root systems 
and elliptic geometry, J. Algebra 43 (1976) 305--327.

\bibitem{ChudnovskySeymour}
M.~Chudnovsky,~P.~Seymour, The structure of claw-free graphs in Surveys in combinatorics 2005, 
London Math. Soc. Lecture Note Ser., Cambridge: Cambridge Univ. Press, 327 (2005), 153--171.


\bibitem{drgsurvey}
E.R.~van~Dam,~J.H.~Koolen,~H.~Tanaka, Distance-regular graphs, 
Electronic J. Combin., Dynamic Survey DS22.

\bibitem{Delsarte}
P.~Delsarte, An algebraic approach to the association schemes of coding theory, 
Philips Res. Rep. Suppl. 10 (1973).

\bibitem{45}
A.L.~Gavrilyuk,~A.A.~Makhnev, Distance-regular graph with the intersection 
array $\{45,30,7;1,2,27\}$ does not exist, Discrete Math. Appl. 23 (2013), 225--244.


\bibitem{godsil-93-paper}
C.D.~Godsil, Geometric distance-regular covers, 
New Zealand J. Math. 22(2) (1993), 31--38.

\bibitem{godsil-93}
C.D.~Godsil, Algebraic combinatorics, Champman and Hall
Mathematics Series, Champman and Hall, New York, 1993.

\bibitem{k6}
A.~Hiraki,~K.~Nomura,~H.~Suzuki, Distance-regular graphs of valency 6 and $a_1=1$, 
J. Alg. Combin. 11(2) (2000), 101--134.

\bibitem{KabanovMakhnev}
V.V.~Kabanov, A.A.~Makhnev, On separated graphs with certain regularity conditions, 
Mat. Sb. 187:10 (1996), 73--86. (In Russian.)

\bibitem{Klin-Pech}
M.~Klin,~C.~Pech, A new construction of antipodal distance-regular
covers of complete graphs through the use of Godsil-Hensel matrices, 
Ars Mathematica Contemporanea, 4 (2011), 205--243.

\bibitem{-m}
J.H.~Koolen,~S.~Bang, On distance-regular graphs with smallest eigenvalue at least $-m$, 
J. Combin. Theory Ser. B 100(6) (2010), 573--584.

\bibitem{koo-park}
J.H.~Koolen,~J.~Park, Distance-regular graphs with $a_1$ or $c_2$ at least half the valency,
J. Combin. Theory Ser. A 119(3) (2012), 546--555.

\bibitem{shilla}
J.H.~Koolen,~J.~Park, Shilla distance-regular graphs, 
European J. Combin. 31(8) (2010), 2064--2073.

\bibitem{kyp}
J.H.~Koolen, J.~Park,~H.~Yu, An inequality involving the second largest and 
smallest eigenvalue of a distance-regular graph, 
Linear Algebra Appl. 434 (2011), 2404--2412.

\bibitem{Metsch}
K.~Metsch, Improvement of Bruck's completion theorem, 
Des. Codes Cryptogr. 1(2) (1991), 99--116.

\bibitem{neumaier-m}
A.~Neumaier, Strongly regular graphs with smallest eigenvalue $-m$,
Arch. Math. (Basel) 33(4) (1979/80), 392--400.

\bibitem{srg25}
A.J.L.~Paulus, Conference matrices and graphs of order $26$, Technische Hogeschool Eindhoven, 
report WSK 73/06, Eindhoven, 1983, 89 pp.

\bibitem{Powers}
D.L.~Powers, Eigenvectors of distance-regular graphs, 
SIAM J. Matrix Anal. Appl. 9 (1988), 399--407.

\bibitem{Spielman} 
D.A.~Spielman, Faster isomorphism testing of strongly regular graphs. 
In Proc. 28th ACM Symp. on Theory of Computing (STOC'96) (1996), 576--584.

\bibitem{Wilmes}
X.~Sun,~J.~Wilmes, Faster canonical forms for primitive coherent configurations. 
In Proc. 47th ACM Symp. on Theory of Computing (STOC'15) (2015), 693--702.

\bibitem{Wilson}
R.M.~Wilson, Nonisomorphic Steiner triple systems, 
Math. Z. 135 (1973/74), 303--313.

\bibitem{yamazaki}
N.~Yamazaki, Distance-regular graphs with $\Gamma_1(x)\simeq 3*K_{a+1}$, 
European J. Combin. 16(5) (1995), 525--536.

\end{thebibliography}
\end{document}